\documentclass[a4paper, 11pt]{amsart}
\usepackage{iftex}
\usepackage{a4wide}
\usepackage[UKenglish]{babel}
\ifPDFTeX 
  \usepackage[utf8]{inputenc}
  \usepackage[T1]{fontenc}
\else 
  \usepackage{fontspec}
\fi
\usepackage{lmodern}
\usepackage{graphicx,subfigure}
\usepackage[colorlinks=true,linkcolor=blue,citecolor=red]{hyperref}
\usepackage{xcolor}
\colorlet{RED}{red}
\usepackage{amssymb, amsfonts, amsmath}
\usepackage[inline]{enumitem}
\usepackage{cases}
\usepackage{lineno}
\usepackage{todonotes}

\usepackage{geometry}
\usepackage{mathtools}
\usepackage{yfonts}
\usepackage[foot]{amsaddr}

\usepackage{geometry}
\usepackage{floatrow}
\usepackage[export]{adjustbox}
\usepackage{comment}
\usepackage{textcomp}
\usepackage{mathrsfs}
\usetikzlibrary{arrows}
\usepackage[ruled,vlined]{algorithm2e}
\usepackage{multirow}
\usepackage{booktabs}
\usepackage{listings}
\usepackage[version=4]{mhchem}
\usepackage{hyperref}
\DeclareCaptionFormat{boldformat}{\textbf{#1#2}#3}

\usepackage{pgfplots}
\pgfplotsset{compat=newest}
\usepackage{tikz}

\theoremstyle{plain}
\newtheorem{thm}{Theorem}[section]

\newtheorem{rem}[thm]{Remark}
 
\theoremstyle{definition} 
\newtheorem{defi}{Definition}[section] 
\theoremstyle{remark}

\newcommand{\R}{\mathbb{R}}
\newcommand{\N}{\mathbb{N}}

\newcommand{\eps}{\varepsilon}

\usepackage{booktabs}
\usepackage{tabularx}

\newtheorem{innercustomthm}{\normalfont \textbf{Theorem}}
\newenvironment{oneshot}[1]
  {\renewcommand\theinnercustomthm{#1}\innercustomthm}
  {\endinnercustomthm}

\definecolor{ffqqqq}{rgb}{1,0,0}
\definecolor{ffzzqq}{rgb}{1,0.6,0}
\definecolor{qqccqq}{rgb}{0,0.6,0}

\allowdisplaybreaks

\makeatletter
\newcommand\thankssymb[1]{\textsuperscript{\@fnsymbol{#1}}}

\author{Tommaso Tenna}
\address[Tommaso Tenna]{Université Côte d’Azur, CNRS, LJAD, Parc Valrose, F-06108 Nice, France}
\email[Tommaso Tenna]{tommaso.tenna@univ-cotedazur.fr}
\date{}

\begin{document}

\title[Projective integration for parabolic systems]{Projective integration schemes for nonlinear degenerate parabolic systems}

\begin{abstract}
A general high-order fully explicit scheme based on projective integration methods is here presented to solve systems of degenerate parabolic equations in general dimensions. The method is based on a BGK approximation of the advection-diffusion equation, where we introduce projective integration method as time integrator to deal with the stiff relaxation term. This approach exploits the clear gap in the eigenvalues spectrum of the kinetic equation, taking into account a sequence of small time steps to damp out the stiff components, followed by an extrapolation step over a large time interval. The time step restriction on the projective step is similar to the CFL condition for advection-diffusion equations. In this paper we discuss the stability and the consistency of the method, presenting some numerical simulations in one and two spatial dimensions.

\end{abstract}


\maketitle


\section{Introduction}
The numerical study of degenerate parabolic equations has gained particular interest in the last decades. These problems arise in several physical applications, such as flow of gas in porous media \cite{diaz1987}, radiative transport \cite{pomraning1979} or cell chemotactic movement \cite{natalini2015}. The numerical approximation of these equations is challenging due to the presence of possible degenerate points, in which the parabolic system loses its properties, exhibiting hyperbolic behavior with finite speed of propagation and sharp fronts. In general, classical solutions might not exist even for $C^2$ initial data and weak solutions must be considered.
Many different schemes have been proposed to solve degenerate parabolic equations, such as local discontinuous Galerkin methods \cite{zhangwu2009}, mixed finite element methods \cite{arbogast1996}, central finite volume schemes \cite{tadmor2000}, finite volume schemes preserving steady-states \cite{bessemoulin2012}, high order finite difference WENO schemes \cite{liushuzhang2011} and WENO schemes based on deep learning techniques \cite{kossaczka2022neural}. Among these methods, a discrete kinetic approximation of these equations have been proposed by Aregba-Driollet, Natalini and Tang \cite{aregba2004}, inspired by relaxation schemes \cite{jinxin1995, aregba2000} and kinetic approximation of hydrodynamic equations in the hyperbolic setting \cite{lionsperthametadmor1994}. The idea relies on a Bhatnagar-Gross-Krook (BGK) approximation of the macroscopic equation, where both the kinetic velocities and the source term depend singularly on the relaxation parameter. Different relaxation systems for the approximation of degenerate parabolic equations have been proposed in \cite{naldipareschi1998, naldipareschitoscani2002}, with high-order extension developed in \cite{cavalli2007}, where the authors coupled ENO and WENO schemes for space discretization with IMEX schemes for time advancement. This approach is inspired by kinetic schemes for the Carleman model, but the numerical integration is directly performed at the macroscopic level. In the context of relaxation schemes, other diffusive kinetic models and approximations have been proposed, see for instance \cite{lionstoscani1997, jinpareschi1998, naldipareschi2000}.\\
The main advantage of relaxation schemes is the possibility to treat the scalar and the system case in the same way at the numerical level, with the feasibility of parallelizing the algorithm due to the diagonal form of the approximating problems. Indeed discrete kinetic approximations are based on linear hyperbolic terms, whose structure is very likely for numerical and theoretical purpose \cite{aregba2000, aregba2004, natalini1996}. The nonlinearity inside the derivatives is replaced by a semilinearity: the differential part becomes linear and all the nonlinearity is concentrated inside the source term. On the other hand, the main disadvantage of this approach is the appearance of a stiff source term, describing the relaxation of the new variable towards the equilibrium, which requires suitable time integrators. \\
In this perspective, the idea is to use an efficient explicit time integrator which is able to speed up the computation due to the particular eigenvalues separation in the spectrum of the kinetic formulation. Projective integration method was introduced in \cite{gear2003} as a technique to accelerate a brute-force integration of stiff systems of ordinary differential equations, which require restrictive conditions on the time step. Such stiff problems are characterized by the presence of different scales, that correspond to a clear gap in the eigenvalue spectrum, in which the fast modes (corresponding to the eigenvalues with large negative part) decay quickly, whereas the slow modes are the active components, still present in the solution. The idea, extended in \cite{lafitte2012} for kinetic equations with a diffusive scaling, is based on an extrapolation in time over a large time step, performed after taking a few inner steps to damp out the transient corresponding to the fast modes. Several extensions have been proposed both for hyperbolic problems \cite{lafittemelis2017} and for high order approximation of kinetic equations \cite{lafittelejon2016}. We mention also telescopic projective integration techniques, where the projective integration method have been extended to deal with problems having multiple time scales, both for ODEs systems \cite{gear2003telescopic} and for kinetic equations \cite{melis2018, melis2019, bailo2022projective}. \\
In \cite{aregba2004, cavalli2007} the methods are based on splitting techniques, restricting the order in time to $2$. The order of convergence in time could be increased by solving explicitly the BGK equation, without projecting over the Maxwellian space. Many asymptotic-preserving schemes based on IMEX techniques have been proposed to integrate stiff kinetic equations and we refer to the recent survey \cite{boscarino2024} for a detailed analysis on the topic.\\
The idea of this paper is to propose a fully explicit and robust projective integration scheme to numerically solve the stiff BGK system with arbitrary order of accuracy in time, extending the idea of \cite{lafittemelis2017} to the case of possibly degenerate parabolic systems.\\
The plan of the paper is the following: in Section \ref{Section::Relaxation} we introduce discrete kinetic relaxation schemes for convection-diffusion equations, investigating different models. We construct the numerical scheme, focusing on the phase space discretization. In Section \ref{Section::Projective} we introduce projective integration method to treat the time discretization with a fully explicit approach, by exploiting the particular structure of the BGK operator. Based on \cite{lafittelejon2016}, we investigate the spectrum of the operator and we perform stability and consistency analysis of the method. In Section \ref{Section::Numerical_Simulations} we present some numerical simulations in both one and two space dimensions for different problems. All the details about the numerical schemes are reported in the Appendix.

\section{Relaxation schemes for nonlinear degenerate parabolic systems}
\label{Section::Relaxation}
In this section, we study discrete kinetic schemes for systems of conservation laws with possibly degenerate diffusion, referring to \cite{aregba2004}, but other kinetic models and approximations can be built, see for instance \cite{jinpareschi1998}. Let $u: \R^D \times \R^+ \to \R^K$, $D >0$, $K>0$, be a weak solution of the following system
\begin{equation}
\label{system_diffusive}
    \partial_t \, u + \displaystyle \sum_{d=1}^D \partial_{x_d} A_{d}(u) = \Delta_x \, B(u), \qquad (x,t) \in \R^D \times (0, +\infty)
\end{equation}
for a given initial condition
\begin{equation}
    u(x,0)=u_0(x), \qquad x \in \R^D.
\end{equation}
Here $t$ defines the time, $x_d$ defines the spatial variable for each dimension $d=1,\dots,D$ and $u_0: \R^D \to \R^K$. The functions $A_d: \R^K \to \R^K$ and $B: \R^K \to \R^K$ are assumed to be Lipschitz-continuous. In addition, for all $u$ lying in a fixed domain $\Omega \in \R^D$ it holds true:
\begin{itemize}
    \item For all $\xi \in \R^D$, $\displaystyle \sum_{d=1}^D \xi_d A'_d(u)$ has real eigenvalues and it is diagonalizable,
    \item The real part of the eigenvalues of $B'(u)$ is non-negative.
\end{itemize}
In \cite{aregba2000} a family of discrete diffusive kinetic schemes has been proposed, in order to construct a numerical approximation of the original system \eqref{system_diffusive}. Let us consider $f^\eps: \R^D \times \R^+ \to \R^L$ as the solution of the following kinetic model in general form:
\begin{equation}
\label{kinetic_approximation}
    \begin{cases}
    \displaystyle \partial_t f^\eps + \sum_{d=1}^D \Gamma_d^\eps \, \partial_{x_d} f^\eps = \frac{1}{\eps} (\mathcal{M}(\eps, u^\eps) - f^\eps),\\
    \displaystyle u^\eps(x,t)= \sum_{l=1}^L f_l^\eps(x,t),
    \end{cases}
\end{equation}
where $L>0$ is the number of discrete velocites and let us take the initial condition
\begin{equation}
    f^\eps(x,0)=f_0^\eps(x).
\end{equation}
Here, $\mathcal{M}: \R^+ \times \R^K \to \R^L$ are the Maxwellians and $\Gamma_d^\eps \in \R^{L \times L}$ are constant diagonal matrices, where the diagonal entries $\gamma_{ld}^\eps$ are given by
\begin{equation}
    \gamma_{ld}^\eps := \lambda_{ld} + \frac{\theta_{ld}}{\sqrt{\eps}},
\end{equation}
for some fixed real constants $\lambda_{ld}$, $\theta_{ld}$. The relaxation parameter $\eps \in \R^+$ can be a physical parameter or an artificial one and the kinetic system \eqref{kinetic_approximation} converges (at least formally) to the original macroscopic system
as $\eps \to 0$, under an analogous of the \textit{subcharacteristic condition}.\\
The Maxwellian $\mathcal{M}$ essentially embeds the macroscopic variable $u^\eps$ in the kinetic space and the right hand side models a relaxation of the kinetic variable $f^\eps$ towards the equilibrium state given by the Maxwellian $\mathcal{M}$. From now on, let us consider $\Omega \in \R^K$ a bounded domain, such that $u(x,t) \in \Omega$ for all $(x,t) \in \R^D \times \R^+$. When $\eps \to 0$, the kinetic approximation formally converges to the macroscopic model under some fundamental assumptions on the Maxwellian functions:
\begin{enumerate}[label=(M{\arabic*})]
\label{Maxwellian_compatibility}
    \item 
    \label{Maxwellian_compatibility_M1}
    $\displaystyle \sum_{l=1}^L \mathcal{M}_l(\eps, w)= w$ for all $\eps \in (0,1]$ and for all $w \in \Omega$,  for all $d= 1,\dots,D$
    \item \label{Maxwellian_compatibility_M2} $\displaystyle \sum_{l=1}^L \gamma_{ld}^\eps \, \mathcal{M}_l(\eps, w)= A_d(w)$ for all $\eps \in (0,1]$ and for all $w \in \Omega$,  for all $d= 1,\dots,D$
    \item 
    \label{Maxwellian_compatibility_M3}
    $\displaystyle \sum_{l=1}^L \theta_{ld} \theta_{lj} \mathcal{M}_l(0, w)= \delta_{dj} B(w)$ for all $w \in \Omega$, for all $d,j = 1,\dots,D$\\
    \item \label{Maxwellian_compatibility_M4} $\mathcal{M}_l(\eps,w) \to \mathcal{M}_l(0,w)$ when $\eps \to 0$ uniformly for $w \in \Omega$.
\end{enumerate}

If $\mathcal{M}$ satisfies the above conditions it is called a local Maxwellian Function for system \eqref{system_diffusive} on $\Omega$. The first property tells us that the projection of the Maxwellian is the identity, which means that considering a given macroscopic variable embedded in the microscopic space and projecting it back, we obtain the original state. The second property and the third property are necessary to guarantee the consistency with the macroscopic fluxes and, therefore, to preserve the structure of the macroscopic model in the limit of the kinetic model \eqref{kinetic_approximation} towards \eqref{system_diffusive}. The last property is fundamental to obtain a uniform convergence with respect to $\eps$, when $\eps \to 0$.\\
The stability condition for this model deals with the following definition.
\begin{defi}
    A local Maxwellian function is a monotone Maxwellian function (MMF) if for all $l \in \{1,\dots,L\}$, for all positive $\eps \leq 1$ and for all $u \in \Omega$, the real parts of the eigenvalues of $\mathcal{M}'_l(u)$ are nonnegative.  
\end{defi}
Rigorous theoretical results about the convergence of the solution to the BGK model \eqref{kinetic_approximation} towards a weak solution to the system \eqref{system_diffusive} have been proven for the scalar case in \cite{bouchut2000} and for a class of one dimensional strongly degenerate parabolic systems in \cite{lattanzionatalini2002}.\\ 
In the next subsection we give a formal derivation of the macroscopic equation, starting from the Chapman-Enskog expansion and assuming the properties \ref{Maxwellian_compatibility_M1}-\ref{Maxwellian_compatibility_M4}.

\subsection{Chapman-Enskog expansion}
Let us define the projection matrix $P \in \R^{L \times K}$ and the auxiliary variables 
\begin{equation}
    u^\eps := P\, f^\eps = \sum_{l=1}^L f_l^\eps, \qquad v_d^\eps = P \, \Gamma_d \, f^\eps = \sum_{l=1}^L \gamma_{ld}^\eps \, f_l^\eps \qquad d=1,\dots,D. 
\end{equation}
Then from \eqref{kinetic_approximation} and the compatibility assumptions \ref{Maxwellian_compatibility_M1}-\ref{Maxwellian_compatibility_M4}, we have
\begin{equation}
\label{system_projected_chapman}
    \begin{cases}
        \partial_t u^\eps +  \displaystyle \sum_{d=1}^D \partial_{x_d} v_d^\eps = 0,\\
        \partial_t v_d^\eps + \displaystyle \sum_{j=1}^D \partial_{x_j} \left(P\, \Gamma_d \Gamma_j f^\eps \right) = \displaystyle \frac{1}{\eps} \left(A_d(u^\eps)-v_d^\eps \right),\\
        \displaystyle \sum_{l=1}^L \theta_{ld} \theta_{lj} f^\eps_l = \delta_{dj} B(u^\eps)
    \end{cases}
\end{equation}
Let us consider a formal Chapman-Enskog expansion of $f^\eps$ in the form 
\begin{equation}
\label{f_approximation_Chapman}
    f^\eps = \mathcal{M}(u^\eps) + \eps \,g^\eps + \mathcal{O} (\eps^2),
\end{equation}
Then, from \eqref{system_projected_chapman}, one has
\begin{equation}
    v_d^\eps= A_d(u^\eps) - \eps \left(\partial_t v_d^\eps + \displaystyle \sum_{j=1}^D \partial_{x_j} \left(P\, \Gamma_d \Gamma_j f^\eps \right) \right) + \mathcal{O}(\eps^2),
\end{equation}
from which we have
\begin{equation}
    v_d^\eps= A_d(u^\eps) - \eps \left( \partial_t v_d^\eps + \displaystyle \sum_{j=1}^D \partial_{x_j} \sum_{l=1}^L \gamma_{dj} \gamma_{lj} f^\eps_l \right) + \mathcal{O}(\eps^2).
\end{equation}
This implies
\begin{equation}
    \partial_t u^\eps + \displaystyle \sum_{d=1}^D \partial_{x_d} A_d(u^\eps) = \eps \sum_{d=1}^D \partial_{x_d} \left( \displaystyle \partial_t v_d^\eps + \displaystyle \sum_{j=1}^D \partial_{x_j} \sum_{l=1}^L \gamma_{dj} \gamma_{dl} f^\eps_l \right) + \mathcal{O}(\eps^2).
\end{equation}
 Then, we have
\begin{equation}
    \partial_t v_d^\eps = A'_d(u^\eps) \, \partial_t u^\eps + \mathcal{O}(\eps) = - \displaystyle \sum_{j=1}^D A'_d(u^\eps)\, A_j'(u^\eps) \partial_{x_j} u^\eps + \mathcal{O}(\eps^2),
\end{equation}
which leads
\begin{multline*}
    \partial_t u^\eps + \displaystyle \sum_{d=1}^D \partial_{x_d} A_d(u^\eps) \\ =\displaystyle \sum_{d=1}^D \displaystyle \sum_{j=1}^D \partial_{x_d}\,\partial_{x_j} \sum_{l=1}^L \theta_{dj} \theta_{dl} f^\eps_l + 2\sqrt{\eps} \displaystyle \sum_{d=1}^D \displaystyle \sum_{j=1}^D \partial_{x_j} \sum_{l=1}^L \lambda_{dj} \theta_{dl} f^\eps_l
    \\+ \eps \left(\displaystyle \sum_{d=1}^D \displaystyle \sum_{j=1}^D \partial_{x_j} \sum_{l=1}^L \lambda_{dj} \lambda_{dl} f^\eps_l  + \sum_{d=1}^D \partial_{x_d} \left( \displaystyle \sum_{j=1}^D A_d'(u^\eps)\,A_j'(u^\eps) \partial_{x_j} u^\eps \right) \right) + \mathcal{O}(\eps^2).
\end{multline*}
Using $(M3)$, when $\eps \to 0$, we finally obtain
\begin{equation}
    \partial_t u^\eps + \displaystyle \sum_{d=1}^D \partial_{x_d} A_d(u^\eps) =\displaystyle \Delta_x B(u^\eps),
\end{equation}
requiring that the leading order terms in the $\eps$-expansion are positive, in the same spirit of the so called \textit{subcharacteristic condition} for the hyperbolic case.
This formal justification shows that the natural convergence rate in \ref{Maxwellian_compatibility_M4} would be $\sqrt{\eps}$. The rigorous proof of this result has been obtained in \cite{bouchut2000}, using a priori bounds on the discrete velocity BGK approximation and kinetic entropy inequalities.

\subsection{Discrete kinetic models} \label{Models_Maxwellians}
In the pure hyperbolic setting, system \eqref{kinetic_approximation} reduces to the one developed in \cite{natalini1996} and \cite{aregba2000}. Let us consider the following hyperbolic system of conservation laws
\begin{equation}
\label{hyperbolic_original}
    \partial_t u + \sum_{d=1}^D \partial_{x_d} A_d(u) = 0.
\end{equation}
and the following approximating discrete BGK model
\begin{equation}
\label{BGK_hyperbolic}
    \begin{cases}
        \partial_t f_l^\eps + \displaystyle \sum_{d=1}^D \lambda_{ld} \partial_{x_d} f_l^\eps = \frac{1}{\eps} ( \tilde{\mathcal{M}}_l(u^\eps)-f_l^\eps), \qquad l=1,\dots,J,\\
        u^\eps = P\,f^\eps = \displaystyle \sum_{l=1}^J f_l^\eps,
    \end{cases}
\end{equation}
where $J$ is the number of velocities needed to approximate the hyperbolic equation, complemented with the compatibility conditions
\begin{equation}
    \sum_{l=1}^J \tilde{\mathcal{M}}_l(u) = u, \qquad \sum_{l=1}^J \lambda_{ld} \tilde{\mathcal{M}}_l(u) = A_d(u), \qquad \forall u \in \Omega, \quad d=1,\dots,D,
\end{equation}
In \cite{natalini1996} it has been shown that the solution of the discrete BGK model \eqref{BGK_hyperbolic}, under the hypothesis of nondecreasing Maxwellian functions, converges to the unique entropy solution of the hyperbolic system \eqref{hyperbolic_original}, see also \cite{natalini1998}.\\
In the same spirit, a model for \eqref{system_diffusive} can be obtained by adding linear combinations of $B$ to the $\tilde{\mathcal{M}}_l$, together with supplementary equations to take into account the diffusive scaling. More precisely, fixing $J' \geq D+1$ the number of velocities you need to deal with the parabolic term and $L=J+J'$, we can design the model as
\begin{equation}
\label{system_discrete_explicit}
    \begin{cases}
        \displaystyle \partial_t f_l^\eps + \sum_{d=1}^D \lambda_{ld}^\eps \partial_{x_d} f_l^\eps = \frac{1}{\eps} (\tilde{\mathcal{M}}_l(u^\eps) - f_l^\eps), \qquad 1 \leq l \leq J,\\[10pt]
        \displaystyle \partial_t f_{J+m}^\eps + \left(\mu + \frac{\theta \sqrt{J'}}{\sqrt{\eps}} \right)\sum_{d=1}^D \sigma_{md}^\eps \partial_{x_d} f_{J+m}^\eps = \frac{1}{\eps} \left(\frac{B(u^\eps)}{J' \theta^2} - f_{J+m}^\eps \right), \qquad 1 \leq m \leq J',\\[10pt]
        \displaystyle u^\eps(x,t)= \sum_{l=1}^L f_l^\eps(x,t),
    \end{cases}
\end{equation}
In particular, for $1 \leq l \leq J$:
\begin{equation}
    \mathcal{M}_l(u) = \tilde{\mathcal{M}}_l(u) + b_l B(u),
\end{equation}
where $b_l \in \R$.
On the other hand, for the last $J'$ equations, we have
\begin{equation}
    \theta \geq 0, \, \mu \geq 0
\end{equation}
and $\{\sigma_d \}_{d=1}^D$ is an orthonormal family in $\{X \in \R^{J'} \text{ s.t. } \displaystyle \sum_{m=1}^{J'} X_m = 0\} = H$.\\
\subsubsection{Diagonal Relaxation Model 1 (DRM1).}
Let us first focus on the one-dimensional case, where we choose $J=2$ and $J'=2$. Then the kinetic model reads as
\begin{equation}
\label{DRM1_kinetic_1D}
    \partial_t f^\eps +  \Gamma^\eps \partial_{x} f^\eps = \frac{1}{\eps} (\mathcal{M}(u^\eps) - f^\eps),
\end{equation}
where 
\begin{equation}
\label{Maxwellians_DRM1}
    \begin{cases}
        \Gamma^\eps=\text{diag} \left(\displaystyle -\lambda, \lambda, -\frac{\mu}{\sqrt{2}}-\frac{\theta}{\sqrt{\eps}}, \frac{\mu}{\sqrt{2}}+\frac{\theta}{\sqrt{\eps}} \right)^T,\\[10pt]
        \mathcal{M}(u)=\displaystyle \frac{1}{2} \left(u-\frac{A(u)}{\lambda}-\frac{B(u)}{\theta^2}, u+\frac{A(u)}{\lambda}-\frac{B(u)}{\theta^2}, \frac{B(u)}{\theta^2}, \frac{B(u)}{\theta^2} \right)^T
    \end{cases}
\end{equation}
Note that the Maxwellian functions do not depend anymore on $\eps$.\\
Let us suppose that for all $u$, $A'(u)$ and $B'(u)$ have a basis of common eigenvectors, denoting $\lambda_k(U)$ and $\theta^2_k(u)$ their respective eigenvalues for $1 \leq k \leq K$ and let us assume that the eigenvalues of $B'$ are real, then $\mathcal{M}$ is a Monotone Maxwellian Function if
\begin{equation}
    \sup_{u \in \Omega} \sup_{1 \leq k \leq K} \frac{|\lambda_k(u)|}{\lambda} + \frac{\theta^2_k(u)}{\theta^2} \leq 1,
\end{equation}
which is a generalized subcharacteristic condition \cite{aregba2004}.

\par
Let us consider now the two-dimensional model, in which $J' \geq D+1=3$, $\theta >0$, $\mu \geq 0$. Let $\{ \sigma^{(1)}, \sigma^{(2)} \}$ be orthonormal in $H$.
\begin{equation}
    \gamma_1^\eps=\begin{pmatrix}
        \lambda_1(-1,0,1)^T\\[7pt]
        \displaystyle \left(\mu+\frac{\theta\sqrt{J'}}{\sqrt{\eps}} \right) \sigma^{(1)}
    \end{pmatrix}, \quad \gamma_2^\eps=\begin{pmatrix}
        \lambda_2(-1,0,1)^T\\[7pt]
        \displaystyle \left(\mu+\frac{\theta\sqrt{J'}}{\sqrt{\eps}} \right) \sigma^{(2)}
    \end{pmatrix}.
\end{equation}
The Maxwellian functions are:
\begin{equation}
    \mathcal{M}(u)=\frac{1}{3} \begin{pmatrix}
        \displaystyle u-2\frac{A_1(u)}{\lambda_1}+\frac{A_2(u)}{\lambda_2}\\[7pt]
        \displaystyle u-2\frac{A_1(u)}{\lambda_1}+\frac{A_2(u)}{\lambda_2}\\[7pt]
        \displaystyle u-2\frac{A_1(u)}{\lambda_1}+\frac{A_2(u)}{\lambda_2}\\[7pt]
        (0,\dots,0)^T
    \end{pmatrix} + \frac{B(u)}{\theta^2} \begin{pmatrix}
        \displaystyle-\frac{1}{3}(1,1,1)^T\\[7pt]
        \displaystyle\frac{1}{J'}(1,\dots,1)^T
    \end{pmatrix}.
\end{equation}
The Maxwellian $\mathcal{M}$ is a MMF in the scalar case ($K=1$) if the following condition is satisfied:
\begin{equation}
    \max_{u \in I} \left( \frac{2\,A_1'}{\lambda_1}-{A_2'}{\lambda_2}, -\frac{A_1'}{\lambda_1}+{2\,A_2'}{\lambda_2}, -\frac{A_1'}{\lambda_1}-{A_2'}{\lambda_2} \right) \leq 1 - \frac{B'}{\theta^2},
\end{equation}
which is equivalent to the subcharacteristic condition presented in \cite{aregba2000}, for the hyperbolic case $\theta=0$.

\subsubsection{Diagonal Relaxation Model 2 (DRM2)}
In this model, we consider the same structure of (DRM1) replacing $\lambda$ by two distinct values $\lambda_m$ and $\lambda_p$. In the one-dimensional case the kinetic model is given by \eqref{DRM1_kinetic_1D}, where now
\begin{equation}
    \begin{cases}
        \Gamma^\eps=\text{diag} \left(\displaystyle \lambda_m, \lambda_p, -\frac{\mu}{\sqrt{2}}-\frac{\theta}{\sqrt{\eps}}, \frac{\mu}{\sqrt{2}}+\frac{\theta}{\sqrt{\eps}} \right)^T,\\[10pt]
        \mathcal{M}_1(u)=\displaystyle \frac{1}{\lambda_p-\lambda_m}\left(\lambda_p \left(u-\frac{B(u)}{\theta^2} \right) -A(u)\right),\\[10pt]
        \mathcal{M}_2(u)=\displaystyle \frac{1}{\lambda_p-\lambda_m}\left(-\lambda_m \left(u-\frac{B(u)}{\theta^2} \right) +A(u)\right), \\[10pt]
        \mathcal{M}_3(u)=\mathcal{M}_4(u) = \displaystyle \frac{B(u)}{2\, \theta^2}.
    \end{cases}
\end{equation}
Supposing w.l.o.g. that $\lambda_m < \lambda_p$, the Maxwellian $\mathcal{M}$ is a MMF if
\begin{equation}
    \lambda_m \left( 1- \frac{B'(u)}{\theta^2} \right) \leq A'(u) \leq \lambda_p \left( 1- \frac{B'(u)}{\theta^2} \right), \qquad \forall \, u \in I.
\end{equation}
This condition allows for a sharper approximation of $A'$ than the previous one, suggesting that the numerical approximation could provide better results \cite{aregba2004}. More specifically, we take
\begin{equation}
    \lambda_m = \inf_{u \in I} \, \frac{A'(u)}{1-B'(u)/\theta^2}, \qquad \lambda_p = \sup_{u \in I} \, \frac{A'(u)}{1-B'(u)/\theta^2},
\end{equation}
to guarantee the maximality of $\lambda_m$ and the minimality of $\lambda_p$.
\subsubsection{A 3 velocities model (OVM)}
For the sake of completeness, we show another choice for the kinetic model, characterized by orthogonal velocities and $\eps$-dependent Maxwellians
\begin{equation}
\label{3velocities_kinetic_1D}
    \partial_t f^\eps +  \Gamma^\eps \partial_{x} f^\eps = \frac{1}{\eps} (\mathcal{M}(u^\eps) - f^\eps),
\end{equation}
where
\begin{equation}
\begin{cases}
    \Gamma^\eps = \text{diag} \left( 0, \lambda, -\lambda \right),\\
    \mathcal{M}_1(\eps, u)= u - \displaystyle \frac{B(u)}{\theta^2},\\[8pt]
    \mathcal{M}_2(\eps, u)= \displaystyle \frac{\left(\lambda + \frac{\theta}{\sqrt{\eps}} \right)}{2\left|\lambda + \frac{\theta}{\sqrt{\eps}} \right|^2} \, A(u) + \frac{B(u)}{\theta^2},\\[8pt]
    \mathcal{M}_3(\eps, u)= -\displaystyle \frac{\left(\lambda + \frac{\theta}{\sqrt{\eps}} \right)}{2\left|\lambda + \frac{\theta}{\sqrt{\eps}} \right|^2} \, A(u) + \frac{B(u)}{\theta^2}.\\[8pt]
\end{cases}
\end{equation}
In this case, if the equation \eqref{hyperbolic_original} does not degenerate or for instance
\begin{equation}
    B'(u) \geq |A'(u)|
\end{equation}
holds, then the monotonicity conditions are verified for $\eps$ small if $\theta^2 > \max B'$. 
\begin{rem}
It is important to observe that by approximating the original system \eqref{system_diffusive} with a kinetic model of the form \eqref{kinetic_approximation}, a modeling error is introduced and it is strictly $\eps$-dependent. However, in any numerical approximation of this model, a finite value of $\eps$ must be taken. This choice is dictated by the consistency of the numerical schemes with respect to the original problem. 
\end{rem}

\subsection{Numerical Approximation}
\label{Subsection_Numerical_Approximation}
The main idea is to avoid splitting techniques to numerically solve \eqref{system_discrete_explicit}. First, we introduce a phase space discretization using a discrete velocity model among the ones introduced above, for a general matrix 
\begin{equation}
    \Gamma^\eps = \text{diag} \displaystyle \left(\lambda_1, \dots, \lambda_J, -\frac{\mu}{\sqrt{2}}-\frac{\theta}{\sqrt{\eps}}, \frac{\mu}{\sqrt{2}}+\frac{\theta}{\sqrt{\eps}}  \right),
\end{equation}
coupled with a finite volume scheme for the advection term. Then, we propose a projective integration scheme for the time discretization.\\
\subsubsection{Velocity discretization}
Let us rewrite explicitly the semi-discrete problem in the one dimensional version:
\begin{equation}
\label{kinetic_system_1D}
    \begin{cases}
        \displaystyle \partial_t f_l^\eps + \lambda_l \partial_x f_l^\eps = \frac{1}{\eps} \left(\mathcal{M}_l(u^\eps) - f_l^\eps \right), \qquad 1 \leq l \leq J,\\[10pt]
        \displaystyle \partial_t f_{J+1}^\eps - \left(\frac{\mu}{\sqrt{2}}+\frac{\theta}{\sqrt{\eps}}\right) \partial_x f_{J+1}^\eps = \frac{1}{\eps} \left( \frac{B(u^\eps)}{2\theta^2} - f^\eps_{J+1}\right),\\[10pt]
        \displaystyle \partial_t f_{J+2}^\eps + \left(\frac{\mu}{\sqrt{2}}+\frac{\theta}{\sqrt{\eps}}\right) \partial_x f_{J+2}^\eps = \frac{1}{\eps} \left( \frac{B(u^\eps)}{2\theta^2} - f^\eps_{J+2}\right),\\[10 pt]
        u^\eps(x,t) =\displaystyle \sum_{l=1}^{J+2} f_l^\eps (x,t).
    \end{cases}
\end{equation}
The Maxwellian function $\mathcal{M}_l$ for the $l$-th equation of the system is chosen as one of the Maxwellian described in Section \ref{Models_Maxwellians}.\\
In a BGK approximation \eqref{kinetic_approximation}, the minimal number $L$ of discrete velocities depends on the spatial dimension of the problem. For the one-dimensional case we consider $L=4$ possible velocity directions, whereas for the two-dimensional case we set $L=6$. It has been shown that increasing the number of discrete velocities $L$ does not provide better results in terms of accuracy \cite{aregba2004} and since the computational complexity for the numerical approximation of the relaxation system clearly depends on the number of velocities, we will take those values $L=4$ (in the $1$D case) and $L=6$ (in the $2$D case) fixed from now on.

\subsubsection{Spatial discretization}
The spatial discretization of Equation \eqref{kinetic_system_1D} is performed by treating the advection term alone. Here we detail the approximation for a second-order (and fourth-order) finite volume scheme, but this can be easily extended to higher order methods. In particular, as we will observe in Section \ref{Subsection_Consistency}, we need a central finite volume scheme for the approximation of the last $J'$ equations.\\ 
Assuming a cartesian grid of the form $[0, L_1] \times [0,L_2]$ for the two-dimensional case, we discretize each dimension into $I_2$ volumes of size $\Delta x^2$. Since the system consists of decoupled equations along each direction, we detail the spatial discretization based on finite volume method in the one-dimensional case, by considering a spatially uniform grid. The domain $[0, L_1]$ is divided in $I$ cells $\mathcal{C}_i = [x_{i-1/2}, x_{i+1/2}]$, where $x_{i \pm 1/2} = x_i \pm \Delta x /2$ with constant size $\Delta x$, over which the cell average of the solution is approximated. 
A numerical flux function $\mathcal{F}$ is defined to approximate the continuous flux at the interface $x_{i + 1/2}$ of each cell, such that the semi-discretized system reads as
\begin{equation}
    \frac{d\,\bar{f}_{l,i}(t)}{dt} = - \frac{1}{\Delta x} \left[ \mathcal{F}_{l, i+1/2}(t) - \mathcal{F}_{l, i-1/2}(t) \right], \qquad i=1,\dots,I \quad l=1,\dots,L.
\end{equation}
This provides an approximation $\bar{f}_{l,i}$ of the cell average $f_{l,i}$, where the numerical flux satisfies for the linear advection with constant velocity $\Gamma_l$
\begin{equation}
    \mathcal{F}_{l, i+1/2}(t) \approx \Gamma_l f_{l, i+1/2}(t), \qquad i=1,\dots,I \quad l=1,\dots,L.
\end{equation}
For the first $J$ equations we consider $k$-th order upwind scheme. The numerical flux for a first order upwind scheme could be defined as
\begin{equation}
    \mathcal{F}_{l, i+1/2}(t) = \begin{cases}
    \Gamma_l \,f_{l, i}\qquad &\text{if } \, \Gamma_l > 0,\\
    \Gamma_l \, f_{l, i+1} \qquad &\text{if } \, \Gamma_l < 0
    \end{cases}
\end{equation}
Higher order extension of upwind schemes are easily recovered (see for instance \cite{shu1998}) and a third order upwind biased scheme reads as
\begin{equation}
    \mathcal{F}_{l, i+1/2}(t) = \begin{cases}
    \displaystyle \Gamma_l \,\frac{2\,f_{l, i+1} + 3\,f_{l, i} - 6\, f_{l, i-1} + f_{l, i-2}}{6} \qquad &\text{if } \, \Gamma_l > 0,\\[10pt]
    \displaystyle \Gamma_l \, \frac{-2\,f_{l, i+2} + 6\,f_{l, i+1} - 3\, f_{l, i} + f_{l, i-1}}{6} \qquad &\text{if } \, \Gamma_l < 0,
    \end{cases}
\end{equation}
for $l=1,\dots,J$. The last $J'$ equations are instead approximated using a centered numerical flux, in order to guarantee the consistency of the scheme with the diffusive BGK scheme \eqref{kinetic_approximation}. The second-order centered flux scheme is given by
\begin{equation}
    \mathcal{F}_{l, i+1/2}(t) = \Gamma_l \, \frac{f_{l, i+1} + f_{l, i-1}}{2}, \qquad l=J+1,\dots,L.
\end{equation}
Fourth order discretization can be obtained by considering the following fluxes 
\begin{equation*}
    \mathcal{F}_{l, i+1/2}(t) = \begin{cases}
    \displaystyle \Gamma_l \,\frac{-f_{l, i+2} + 8\,f_{l, i+1} - 8\, f_{l, i-1} + f_{l, i-2}}{12} \qquad &\text{if } \, \Gamma_l > 0,\\[10pt]
    \displaystyle \Gamma_l \, \frac{f_{l, i+3} - 6\,f_{l, i+2} + 18 f_{l,i+1}- 10\, f_{l, i} - 3\,f_{l, i-1}}{12} \qquad &\text{if } \, \Gamma_l < 0,
    \end{cases} \qquad l=1,\dots,J.
\end{equation*}
\begin{equation*}
    \mathcal{F}_{l, i+1/2}(t) = \Gamma_l \, \frac{f_{l, i-2} - 8\, f_{l, i-1} + 8\,f_{l, i+1} - f_{l, i+2}}{12} \qquad l=J+1,\dots,L.
\end{equation*}
\begin{rem}
    The violation of a maximum principle for a centered flux scheme could lead to unphysical oscillations, resulting in an unstable approximation of the advection equation. However, the consistency in $\eps$ of the kinetic approximation \eqref{kinetic_approximation} with the advection-diffusion equation \eqref{system_diffusive} quickly stabilize the oscillations as $\eps \to 0$, as already observed in \cite{lafitte2012}.
\end{rem}

\section{Projective Integration}
\label{Section::Projective}
\begin{figure}
\begin{tikzpicture}[line cap=round,line join=round,>=triangle 45,x=1cm,y=1cm]
\clip(-1.117015993609361,-0.9149876034609) rectangle (15.454362892408566,5.573921109952696);
\draw [->,line width=1.2pt] (0,0) -- (0,5.14);
\draw [->,line width=1.2pt] (0,0) -- (13.14,0);
\draw [line width=1.1pt,dash pattern=on 1pt off 1pt] (2.78,2.8)-- (6.5,3.66);
\draw [shift={(5.67747899159664,-12.097310924369756)},line width=1.1pt]  plot[domain=1.3393525936765465:1.8978609152431969,variable=\t]({1*15.182123393302678*cos(\t r)+0*15.182123393302678*sin(\t r)},{0*15.182123393302678*cos(\t r)+1*15.182123393302678*sin(\t r)});
\draw [line width=1.1pt,dotted] (2,0)-- (2,4.2);
\draw [line width=1.1pt,dotted] (6.5,0)-- (6.5,3.66);
\draw [line width=1.1pt] (9.16,2.68)-- (12,2);
\draw [line width=1.1pt,dash pattern=on 1pt off 1pt] (7.28,3)-- (11,2.44);
\draw[line width=1.1pt,dotted] (11,0)-- (11,2.44);
\draw[line width=1.1pt] (0,-1)-- (0,0);
\draw[line width=1.1pt] (0,0)-- (-1,0);
\draw (12.4,0.9) node[anchor=north west] {\Large time};
\draw (12,2.2) node[anchor=north west] {\Large $f(t)$};
\begin{scriptsize}
\draw [color=black] (2,0)-- ++(-2pt,0 pt) -- ++(4pt,0 pt) ++(-2pt,-2pt) -- ++(0 pt,4pt);
\draw[color=black] (2,-0.35) node {\Large $t^{n-1}$};
\draw[color=black] (2,4.7) node {\Large $f^{n-1}$};
\draw[color=black] (3.35,2.4) node {\large $f^{n-1, K}$};
\draw [color=black] (6.5,0)-- ++(-2pt,0 pt) -- ++(4pt,0 pt) ++(-2pt,-2pt) -- ++(0 pt,4pt);
\draw[color=black] (6.5,-0.35) node {\Large $t^n$};
\draw[color=black] (6.5,4.1) node {\Large $f^n$};
\draw [color=black] (11,0)-- ++(-2pt,0 pt) -- ++(4pt,0 pt) ++(-2pt,-2pt) -- ++(0 pt,4pt);
\draw[color=black] (11,-0.35) node {\Large $t^{n+1}$};
\draw[color=black] (11,3.1) node {\Large $f^{n+1}$};
\draw [fill=black] (2,4.2) circle (1.8pt);
\draw [fill=black] (2.4,3.04) circle (1.8pt);
\draw [fill=black] (2.78,2.8) circle (1.8pt);
\draw [fill=black] (3.2,2.88) circle (1.8pt);
\draw [fill=black] (6.5,3.66) circle (1.8pt);
\draw [fill=black] (6.88,3.24) circle (1.8pt);
\draw [fill=black] (7.28,3) circle (1.8pt);
\draw [fill=black] (7.68,2.95) circle (1.8pt);
\draw [fill=black] (11,2.44) circle (1.8pt);
\end{scriptsize}
\end{tikzpicture}
\label{PI_Sketch}
\caption{Idea of Projective Integration method.}
\end{figure}
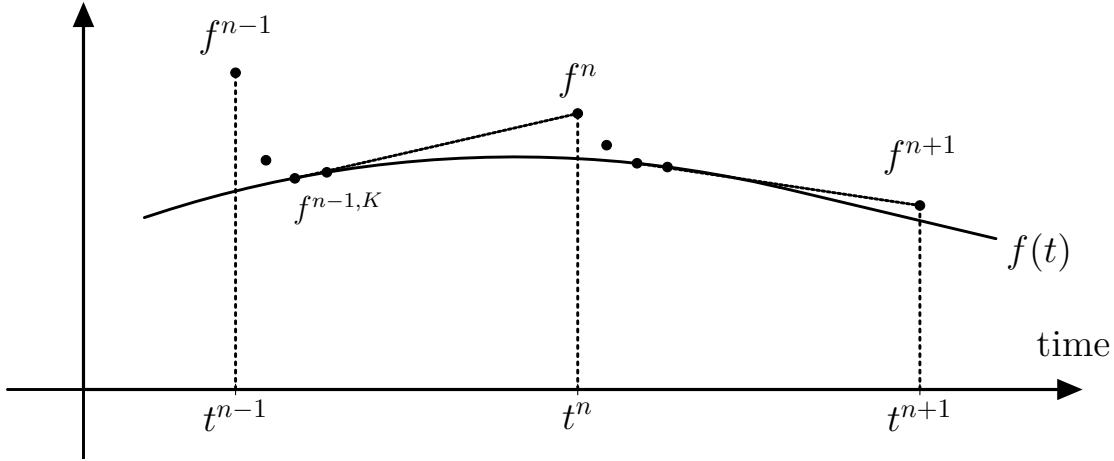

The main idea behind the construction of projective integration methods is to combine a few small time steps of a naive (\textit{inner}) time-stepping method, with a projective step (\textit{outer}), in which an extrapolation over a large time interval is performed. The idea is sketched in Fig. \ref{PI_Sketch}.\\
In \cite{lafitte2012}, the authors show how the choice of an outer time step much larger than the inner time step leads to a splitting of the stability region of the method in two circle-like connected components. In that case, the projective integration parameters can be taken in order to have all the fast eigenvalues in the first stability region and all the slow (dominant) eigenvalues in the dominant stability region.\\
The natural choice for the inner step is to use a simple, explicit method, such as a direct Euler discretization. This approach allows to obtain a fully explicit, arbitrary order in time method for stiff systems, sharing important features with asymptotic-preserving schemes. Indeed, despite the impossibility of evaluating the scheme for $\eps=0$, the computational cost of this method is often independent from the stiffness of the problem. This implies that it is possible to solve stiff problems without requiring restrictive CFL condition, as done in the series of papers \cite{bailo2022projective, melis2019,melis2018} for different stiff kinetic equations.\\
The idea is to use the same approach to numerically approximate the BGK approximation \eqref{kinetic_approximation} in the limit for $\eps \to 0$ with an asymptotic-preserving scheme.\\
Let us consider the semidiscrete version of the system \eqref{kinetic_approximation} obtained through the phase space discretization of the form:
\begin{equation}
\label{semidiscrete}
	\frac{d\,f}{dt} = D_t^\eps (f).
\end{equation}
In our setting, $D_t^\eps$ is the semidiscrete operator given by
\begin{equation}
	\label{semidiscrete_operator}
	D_t^\eps(f)=-\Phi_{\alpha}(f) + \frac{1}{\eps} \left( \mathcal{M}(u)-f \right),
\end{equation}
where $\Phi_\alpha$ is a suitable discretization of the convective derivative, as discussed in Section \ref{Subsection_Numerical_Approximation}.\\
Let us now define step sizes $\delta t$ and $\Delta t$, which represent the inner and the outer time steps respectively. Then, the solution $f(t)$ is approximated by $f^{n,k}$ at the time $t=n\Delta t + k \delta t$.\\
As inner integrator, we choose an explicit time integration scheme $S_{\delta t}$. In particular, the inner step could be computed through the forward Euler method, obtaining
\begin{equation}
	\label{inner_integrator_projective}
	f^{n,k+1}=S_{\delta t}(f^{n,k})= f^{n,k} + \delta t D_t^\eps ( f^{n,k}) \qquad k=0,1,\dots.
\end{equation}
The inner step has to be chosen as $\delta t= \mathcal{O}(\eps)$ to ensure the stability of the method. This is problematic for $\eps \to 0$, since the condition becomes too restrictive and the spatial discretization in the kinetic space may lead unstabilities in the macroscopic space.\\
When $\eps \to 0$, the system formally converges to a limiting equation for which every finite volume scheme needs to sastisfy a weaker stability restriction. The idea of projective integration schemes is to exploit the structure of the kinetic operator to accelerate the integration.\\
Starting from a computed numerical solution $f^n$ at time $t^n=n\Delta t$, one first takes $K+1$ inner steps of size $\delta t$, using the inner integrator \eqref{inner_integrator_projective} and then approximate the time derivative of $f$, computing $f^{n+1}$ via extrapolation in time. The last two updates of \eqref{inner_integrator_projective}, $f^{n,K}$ and $f^{n,K+1}$ are used to project the solution forward
\begin{equation}
    \label{PFE}
    \begin{aligned}
    f^{n+1}=f^{n,K+1}+(\Delta t - (K+1)\delta t) \displaystyle \frac{f^{n,K+1}-f^{n,K}}{\delta t} &\\ =f^{n,K+1}+ M\delta t \displaystyle \frac{f^{n,K+1}-f^{n,K}}{\delta t},&
    \end{aligned}
\end{equation}
in which we have defined the relative size of the extrapolation $M$ as:
\begin{equation}
	M=\frac{\Delta t}{\delta t}-(K+1).
\end{equation}
This method \eqref{PFE} is called \textit{projective forward Euler} (PFE), due to the fact that the extrapolation is performed using a first-order approximation of the derivative.\\ 
\subsubsection*{High-order Projective Runge-Kutta methods.} To obtain higher-order accuracy in time we can employ any Runge-Kutta method as inner integrator, replacing each time derivative evaluation in a classical Runge-Kutta method by $K+1$ steps of an inner integrator, as discussed in \cite{lafittelejon2016}.\\
However, the order of the scheme is dominated by the outer integrator. The idea is to use a particular higher order extension of projective integration, based on Runge-Kutta methods. Here, we briefly recall the construction of the projective Runge-Kutta (PRK) method. Using the same coefficients introduced in Appendix \ref{Appendix_RK} for classical Runge-Kutta methods, for an explicit S-stage PRK scheme, we compute the values $\mathbf{k}_s$ as  \\
\begin{align}
\label{PRK_scheme}
    s = 1 : &\begin{cases}
        f^{n, k+1} = f^{n, k} + \delta t \, D_t(f^{n, k} ) \qquad 0 \leq k \leq K,\\[10pt]
        \mathbf{k}_1 = \displaystyle \frac{f^{n,K+1}-f^{n,K}}{\delta t},
    \end{cases}\\[10pt]
\label{PRK_scheme_2}
    2 \leq s \leq S : &\begin{cases}
        f^{n +c_s, 0} = \displaystyle f^{n, K+1} + \left(c_s \Delta t - (K+1)\delta t \right) \, \sum_{l=1}^{s-1} \frac{a_{s,l}}{c_s} \, \mathbf{k}_l,\\[10pt]
        f^{n+c_s, k+1} = \displaystyle f^{n+c_s, k} + \delta t \, D_t(f^{n+c_s, k} ) \qquad 0 \leq k \leq K,\\[10pt]
        \mathbf{k}_s = \displaystyle \frac{f^{n+c_s,K+1}-f^{n+c_s,K}}{\delta t},
    \end{cases}
\end{align}
where $f^{n+c_s, k} \approx f(t^n + c_s\Delta t + k\delta t)$, with $a_{s,l}$ and $c_s$ the \textit{Runge-Kutta matrix} and the \textit{Runge-Kutta nodes}, respectively.
Finally, the numerical solution at time $t^{n+1}$ is computed as
\begin{equation}
    f^{n+1} = f^{n, K+1} + \left( \Delta t - (K+1) \delta t \right) \displaystyle \sum_{s=1}^S b_s \, \mathbf{k}_s,
\end{equation}
where $b_s$ are the \textit{Runge-Kutta weights}. The explicit expression for the coefficients of the Runge-Kutta schemes is given in the Appendix \ref{Appendix_RK}. In Fig. \ref{PRK_scheme_figure} we sketch the idea of Projective Runge-Kutta schemes for a second order case.
\begin{figure}
\begin{tikzpicture}[line cap=round,line join=round,>=triangle 45,x=1cm,y=0.8cm, scale=1.5]
	\clip(-1.5,-1) rectangle (15.454362892408566,5.573921109952696);
	\draw [->,line width=1.2pt] (0,-0.3807087091449403) -- (-0.012971764851690253,4.206590254843377);
	\draw [->,line width=1.2pt] (-0.42618873147435454,0) -- (6.596795947812606,0);
	\draw [shift={(5.196506491400063,5.080337715759148)},line width=1.2pt]
	
	plot[domain=3.4079859364131755:4.884962253444296,variable=\t]({1*4.284630026256087*cos(\t r)+0*4.284630026256087*sin(\t r)},{0*4.284630026256087*cos(\t r)+1*4.284630026256087*sin(\t r)});
	
	\draw [line width=1.2pt] (1.7567814683606504,2.5256662637380183)-- (2.753202833755961,1.248188226681539);
	\draw [line width=1.4pt,dash pattern=on 1pt off 1pt,color=ffqqqq] (1.3704104177110121,2.5817838062369645)-- (1.8996976756144677,1.77392430733169);
	\draw [line width=1.4pt,dash pattern=on 1pt off 1pt,color=ffqqqq] (2.91,1.75)-- (3.96,1.2);
	\draw [line width=1.2pt] (2.753202833755961,1.248188226681539)-- (4,0.5);
	\draw [line width=1.2pt,dotted] (2.753202833755961,1.248188226681539)-- (4.008074863290168,-0.3058301605557321);
	\draw [line width=1.2pt,dotted] (1.5415758353305784,3.3628975058880344)-- (1.5415758353305784,0);
	\draw(6.326354936891005,0.565114976335994) node[anchor=north west] {\Large $\text{time}$};
	\begin{scriptsize}
		\draw [fill=qqccqq] (1.5415758353305784,3.3628975058880344) circle (1.2pt);
		\draw[color=qqccqq] (1.4958977532154432,3.7226124025447254) node {\large $f^n$};
		\draw [fill=qqccqq] (1.6665669562098817,2.9792244502142444) circle (1.2pt);
		\draw [fill=qqccqq] (1.7567814683606504,2.5256662637380183) circle (1.2pt);
		\draw [fill=qqccqq] (1.9288955395354348,2.3089169376834744) circle (1.2pt);
		\draw[color=qqccqq] (2.4551374776332855,2.637757952310261) node {\large $f^{n,K+1}$};
		\draw [fill=ffzzqq] (2.774884052439232,1.2388666875342405) circle (1.2pt);
		\draw [fill=ffzzqq] (3.4, 1.19) circle (1.2pt);
		\draw[color=ffzzqq] (2.2579131623923397,1.1701911736152546) node {\large $f^{n+c_s,0}$};
		\draw [fill=ffzzqq] (2.980435421957341,1.3073838107069433) circle (1.2pt);
		\draw [fill=ffzzqq] (3.213555943366573,1.2821853246683275) circle (1.2pt);
		\draw [fill=black] (4,0.5) circle (1.2pt);
		\draw[color=black] (4.512983363483837,0.5421175445321434) node {\large $  f^{n+1}_{\text{PRK2}}$};
		\draw [fill=black] (4,-0.306) circle (1.2pt);
		\draw[color=black] (4,-0.65) node {\large $ f^{n+1}_{\text{PFE}}$} ;
		\draw [color=black] (4,0)-- ++(-1.5pt,0 pt) -- ++(3pt,0 pt) ++(-1.5pt,-1.5pt) -- ++(0 pt,3pt);
		\draw[color=black] (4.156646036422076,0.2249768232322224) node { $t^{n+1}$};
		\draw [color=black] (1.5415758353305784,0)-- ++(-1.5pt,0 pt) -- ++(3pt,0 pt) ++(-1.5pt,-1.5pt) -- ++(0 pt,3pt);
		\draw[color=black] (1.4958977532154432,-0.26249768232322224) node {\large $t^n$};
	\end{scriptsize}
	
\end{tikzpicture}
\label{PRK_scheme_figure}
\caption{Idea of the Projective Runge-Kutta Integration method. The scheme first computes $f^{n,K+1}$}
\end{figure}
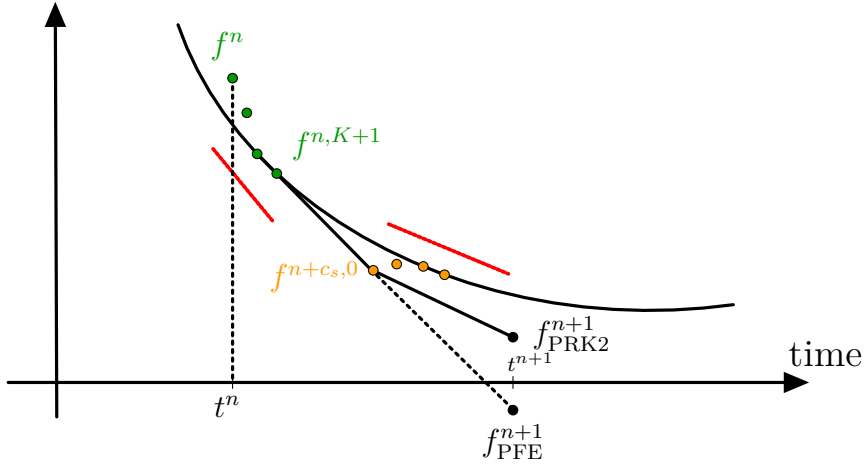

\subsection{Stability analysis}
\label{spectrum_inner_integrator}
The main advantage of projective integration approach is its great speed up for problems with great scale separation, but it is fundamental to estimate the spectrum of the inner integrator for the equation to guarantee the stability of the scheme. The estimation of the spectrum for a general BGK operator in the discrete kinetic approximation framework could be extremely hard. For this reason, we propose here an analysis of the spectrum for a simplified setting, namely taking $A(u)=u$ and $B(u)=u$ in \eqref{system_diffusive} and the diagonal relaxation model (DRM1).\\
Let us consider the scalar BGK approximation in Equation \eqref{DRM1_kinetic_1D} with a linear Maxwellian in the form of \eqref{Maxwellians_DRM1}
\begin{equation}
    \mathcal{M}_{\lambda,\theta} (u)= \left(u + \frac{u}{\lambda} - \frac{u}{\theta^2}, \, u - \frac{u}{\lambda} - \frac{u}{\theta^2}, \, \frac{u}{\theta^2}, \, \frac{u}{\theta^2} \right).
\end{equation}
We can rewrite the semidiscrete kinetic equation \eqref{semidiscrete_operator} in the Fourier (spatial) domain as
\begin{equation}
    \partial_t \hat{f}(\xi_k) = \mathcal{K} \hat{f}(\xi_k), \qquad \mathcal{K} = \frac{1}{\eps} \left( -\eps D + MP - I  \right)
\end{equation}
where $\hat{f} \in \R^4$, $I \in \R^{4\times4}$ is the identity matrix, $P \in \R^{4\times4}$ is the orthogonal projection matrix and $M \in \R^{4 \times 4}$ is given by
\begin{equation}
M = \frac{1}{2}\begin{pmatrix}
        1+\displaystyle \frac{1}{\lambda}-\displaystyle \frac{1}{\theta^2} & 0 & 0 & 0 \\
        0 & 1-\displaystyle \frac{1}{\lambda}- \displaystyle\frac{1}{\theta^2} & 0 & 0 \\
        0 & 0 & \displaystyle\frac{1}{\theta^2} & 0 \\
        0 & 0 & 0 & \displaystyle\frac{1}{\theta^2}
    \end{pmatrix}.
\end{equation}
The advection matrix $D \in \R^{4 \times 4}$ can be written as
\begin{equation}
\label{D_matrix_spectrum}
    D = \begin{pmatrix}
        \alpha+i\,\beta & 0 & 0 & 0 \\
        0 & \alpha-i\,\beta & 0 & 0 \\
        0 & 0 & \frac{1}{\sqrt{\eps}}\left(\xi+i\,\gamma\right) & 0 \\
        0 & 0 & 0 & \frac{1}{\sqrt{\eps}}\left(\xi-i\,\gamma\right)
    \end{pmatrix}.
\end{equation}
Here $D$ represents the Fourier matrix of the spatial discretisation chosen for the advection term, and $\alpha$, $\beta$ are respectively the real and the imaginary part of the upwind approximation
\begin{equation}
\label{Upwind_Neumann}
    \begin{cases}
        \alpha = -\displaystyle\frac{|\lambda|}{6\,\Delta x}\left(3-4\cos(\zeta) + \cos(2\zeta)\right),\\[10pt]
        \beta=-\displaystyle\frac{\lambda}{6\,\Delta x} (8\,\sin(\zeta)-\sin(2\,\zeta)),
    \end{cases}
\end{equation}
and $\xi$ and $\gamma$ are respectively the real and the imaginary part of the centered approximation
\begin{equation}
\label{Centered_Neumann}
    \begin{cases}
        \xi = 0,\\[10pt]
        \gamma=\displaystyle\left(\frac{\mu}{\sqrt{2}}+\frac{\theta}{\sqrt{\eps}} \right) \, \frac{8\sin(\zeta)-\sin(2\,\zeta)}{6\,\Delta x}.
    \end{cases}
\end{equation}

Let us define an auxiliary matrix $\mathcal{A}=-\eps\,D+MP$, such that $\mathcal{K}=\eps^{-1}(\mathcal{A}-I)$.
Since the matrix has the following structure
\begin{equation}
    \mathcal{A}=D+\,u\,e^T,
\end{equation}
where $e$ is the canonical vector in $\R^4$, we can use Sherman-Morrison formula to compute the determinant $\chi_{\mathcal{A}}$ \cite{golubvanloan2007}, which is given by
\begin{equation}
    \chi_\mathcal{A}(\zeta) = \prod_{i=1}^4 (-\eps D_i - \zeta) \, \left(1+\frac{1}{4}\,\sum_{j=1}^4 \frac{M_j}{-\eps D_j - \zeta} \right),
\end{equation}
where $D_i$ denotes the $i$-th diagonal entries of $D$ and $M_j$ denotes the $j$-th diagonal entries of $M$.\\
Following \cite{lafittelejon2016}, it is possible to localize the eigenvalues, by using Rouché's Theorem. The main differences here are given by the non-constant $\eps$-dependency of the coefficients $D_j$ and the non-uniform structure of $M_j$. 

\begin{thm}
\label{thm:Spectrum_Inner}
The spectrum of the matrix $\mathcal{A}$ satisfies
\begin{equation}
    \text{Sp}(\mathcal{A}) \subset \mathcal{D} \left(0, C\sqrt{\eps}(|\xi|+|\gamma|) \right) \cup \left\{ \zeta(\eps)\right\},
\end{equation}
where $C$ is a constant depending on the entries of $D$, $\xi$ and $\gamma$ are defined in \eqref{D_matrix_spectrum} and $\zeta(\eps)$ is the dominant eigenvalue.
\end{thm}

For the sake of completeness we retrace the proof of Theorem \ref{thm:Spectrum_Inner} in the Appendix \ref{Appendix_Stability}, extending the results in \cite{lafittelejon2016}. As a consequence, the spectrum of $\mathcal{K}$ satisfies
\begin{equation}
\label{Region_Spectrum}
    \text{Sp} \left(\mathcal{K}\right) \subset \left(\mathcal{D}\left(-\frac{1}{\eps}, C\,\frac{1}{\sqrt{\eps}} \left(|\xi|+|\gamma|\right)\right) \right) \cup \left\{\zeta(\eps)\right\},
\end{equation}
where $\zeta$ is the dominant eigenvalue, whose explicit expression, in the case of an upwind scheme \eqref{Upwind_Neumann} for the hyperbolic part and a centered scheme \eqref{Centered_Neumann} for the parabolic part, is given by
\begin{equation}
\begin{cases}
    \text{Re} \displaystyle \left(\zeta(\eps)\right) = \left(1-\frac{1}{\theta^2} \right) \, \alpha + \mathcal{O} (\eps), \\[10pt]
    \text{Im}\left(\zeta(\eps)\right) = \displaystyle \frac{\beta}{\lambda} + \mathcal{O}(\eps).
\end{cases}   
\end{equation}

We first observe that the spectrum of $\mathcal{K} = \eps^{-1} \left(\mathcal{A} - I \right)$ consists of two well separated clusters for $\eps \to 0$ and its localization on the real axis is given by two contributions: the dominant eigenvalue, whose real part is $\mathcal{O}(1)$, and the remaining eigenvalues, contained in a region around $-1/\eps$. This is coherent with the derivation of the spectrum eigenvalues for a linearized BGK operator, as shown in \cite{cercignani1988}. 

\begin{figure}[h!]
    \begin{center}
    \includegraphics[width=0.95\linewidth]{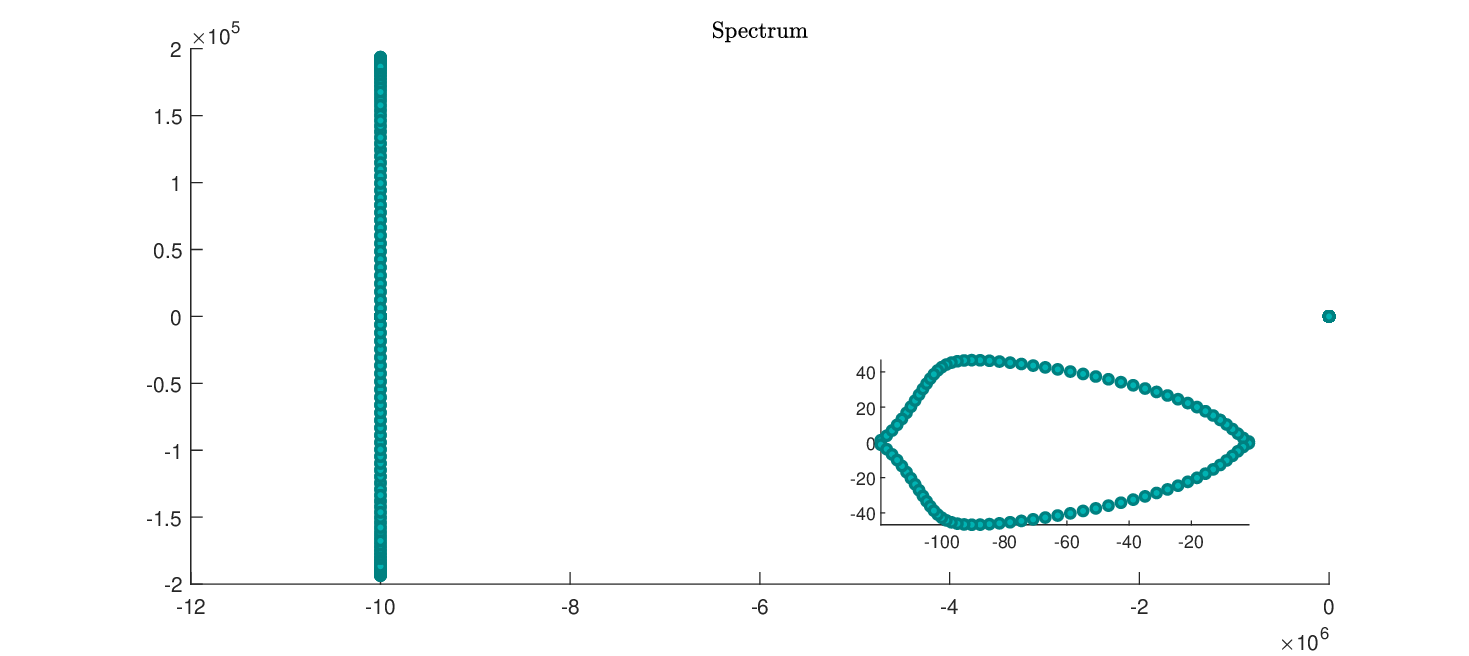}  
    \end{center}
    \caption{Spectrum of the operator for $\eps=10^{-7}$.}
    \label{Spectrum_1e-7}
\end{figure}

Moreover, we assess our theoretical result by introducing a numerical approximation of the spectrum, as done in \cite{bailo2022projective}. The approach is based on a discretization of the problem in phase space, obtaining a semi-discrete system \eqref{semidiscrete}. The set of eigenvalues of the Jacobian of $D_t^\eps$ is precisely the spectrum of our problem.\\
Introducing a finite difference approximation and implementing, for instance, the forward Euler time discretization of \eqref{semidiscrete}, the semidiscrete operator could be rewritten as
\begin{equation}
	G(f^{n,k})=\frac{F(f^{n,k})-f^{n,k}}{\delta t},
\end{equation}
where $F(f^{n,k})=f^{n,k+1}$.
The function $G(f^{n,k})$ represents the discretization of the temporal derivative, which, according to \eqref{semidiscrete}, gives us a discrete of $D_t^\eps$. The Jacobian of this operator and its eigenvalues can be now computed numerically, to obtain a numerical approximation of the spectrum of our problem.\\
In Fig. \ref{Spectrum_1e-7} we show the spectrum of the operator for $\eps=10^{-7}$ with a central difference scheme for the parabolic terms. Similarly to the spectrum of a BGK equation in the diffusive scaling, here the imaginary part is governed by $\sqrt{\eps}$, in accordance to \eqref{Region_Spectrum}.

In \cite{gear2003} it has been shown that the stability region of the PFE consists of the union of two separated disks $\mathcal{D}_1^{\text{PFE}} \cup \mathcal{D}_2^{\text{PFE}}$, where
\begin{equation}
    \mathcal{D}_1^{\text{PFE}}= \mathcal{D} \left(1-\frac{\delta t}{\Delta t}, \frac{\delta t}{\Delta t} \right), \qquad \mathcal{D}_2^{\text{PFE}} \left(0, \left(\frac{\delta t}{\Delta t} \right)^{\frac{1}{K}} \right),
\end{equation}
and $\mathcal{D}(c, r)$ denotes the ball with center $(c,0)$ and radius $r$ in the complex $\lambda$-plane. The eigenvalues lying in $\mathcal{D}_1^{\text{PFE}}$ corresponds to the slowly decaying modes, whereas the eigenvalues $\mathcal{D}_2^{\text{PFE}}$ corresponds to the fast modes, quickly damped by the inner integrator. Since the dominant set of eigenvalues is given by $\mathcal{D}_1^{\text{PFE}}$, the approach of projective integration schemes is to suitably choose the method parameters $\delta t$, $\Delta t$ and $K$ to allow for accurate integration of the modes in $\mathcal{D}_1^{\text{PFE}}$, while ensuring stability for the modes in $\mathcal{D}_2^{\text{PFE}}$.\\
In \cite{lafittelejon2016}, it has been shown that the choice of stable parameters for the PFE guarantees also the stability of high order projective integration methods. Thus, according to \ref{thm:Spectrum_Inner}, we set the parameters for projective integration method as follows:
\begin{equation}
    \delta t = \eps, \qquad K \geq 2.
\end{equation}
The $\Delta t$ is chosen according to the parabolic CFL condition.\\
We refer to the Appendix \ref{stability_PRK_appendix} for a detailed description of the stability analysis of PRK schemes.

\subsection{Consistency analysis}
\label{Subsection_Consistency}
The consistency analysis for our numerical method is performed again in the scalar one-dimensional setting, following the derivation in \cite{lafitte2012}. For this purpose, we consider the expression \eqref{f_approximation_Chapman}. The scheme is given by
\begin{equation}
    f^{n+1}=S_{\delta t} f^n = f^n + \delta t \left(-\Phi^\eps(f^n)+ \frac{\mathcal{M}(u^n)-f^n}{\eps} \right),
\end{equation}
where $\Phi^\eps$ is the spatial discretization operator.\\
If we consider $\delta t = \eps$ (according to the stability condition) and a centered flux for the last $J'$ equations, the inner scheme is then reduced to the following expression
\begin{equation*}
    f^{n+1} = S_{\eps} f^n = \mathcal{M}(u^n) - \eps \Phi^\eps(f^n).
\end{equation*}
The corresponding density $u^{N+1}$ at time $t^{N+1}$ satisfies 
\begin{equation}
    u^{N+1} = u^{N,K+1} + \left(\Delta t - (K+1)\eps \right) \frac{u^{N,K+1}-u^{N,K}}{\eps},
\end{equation}
where $K+1$ is the number of inner time steps taken to guarantee the stability of the method.\\
The truncation error could be defined for the macroscopic quantities as 
\begin{equation}                
    E^{N+1}=\frac{\tilde{u}^{N+1}-u^{N+1}}{\Delta t},
\end{equation}
where $\tilde{u}^N$ is the exact solution for $u$ at time $t^N$.\\
By considering the projective step, we can rewrite the truncation error as 
\begin{equation}
    E^{N+1} = \frac{\tilde{u}^{N+1} - u^{N,K+1}}{\Delta t} - \left(\frac{\Delta t - (K+1)\eps}{\Delta t} \right) \frac{u^{N,K+1}-u^{N,K}}{\eps}.
\end{equation}
This means that the local truncation error of the projective scheme depends on its inner integrator too. Thus, we introduce the local truncation error of the inner integrator related to the kinetic variable $f$ as
\begin{equation}
    e^{N,k+1}_f = \frac{\tilde{f}^{N,k+1}-f^{N,k+1}}{\delta t} = \frac{\tilde{f}^{N,k+1}-f^{N,k+1}}{\eps}.
\end{equation}
We can rewrite both quantities in terms of their solution at time $t^{N,k}$. For $f^{N,k+1}$ we simply use the expression of the inner integrator, which, in this case, is chosen as the Forward Euler scheme
\begin{equation}
    f^{N,k+1} = S_{\delta t} f^{N,k} = \mathcal{M}(u^{N,k}) - \eps \Phi^\eps(f^{N,k}).
\end{equation}
The exact solution $\tilde{f}^{N,k+1}$ is expanded through a Taylor series around $t^{N,k}$, using the kinetic formulation of the model
\begin{equation}
    \begin{aligned}
    \tilde{f}^{N,k+1} &= \tilde{f}^{N,k} + \eps \partial_t \tilde{f}^{N,k} + \mathcal{O}(\eps^2) \\
    &=\tilde{f}^{N,k} - \eps \left( \sum_{d=1}^D \Gamma_d \partial_{x_d} \tilde{f}^{N,k} \right) + \left(M(\tilde{u}^{N,k})-\tilde{f}^{N,k} \right) + \mathcal{O}(\eps^2) \\
    &=S_{\delta t} (\tilde{f}^{N,k}) + \eps \left(\Phi(\tilde{f}^{N,k}) -  \sum_{d=1}^D \Gamma_d \partial_{x_d} \tilde{f}^{N,k} \right) + \mathcal{O}(\eps^2).
    \end{aligned}
\end{equation}
The third equality is obtained from adding and subtracting $\Phi(\tilde{f}^{N,k})$ to the expression and using the forward Euler operator $S_{\delta t}$.\\
This leads to 
\begin{equation}
    e^{N,k+1}_f = S_{\delta t} (e^{N,k}_f) + \left(\Phi(\tilde{f}^{N,k}) - \sum_{d=1}^D \Gamma_d \partial_{x_d} \tilde{f}^{N,k} \right) + \mathcal{O}(\eps).
\end{equation}
By recursion we could rewrite 
\begin{equation}
    e^{N,k+1}_f = \sum_{k=0}^K S^{k}_{\delta t} \left(\Phi(\tilde{f}^{N,K-k}) - \sum_{d=1}^D \Gamma_d \partial_{x_d} \tilde{f}^{N,K-k} \right) + \mathcal{O}((K+1)\eps).
\end{equation}
Let us split the expression as
\begin{multline}
    \label{err_step1}
    e^{N,k+1}_f = \sum_{k=0}^K S^{k}_{\delta t} \left(\mathcal{A}(\tilde{f}^{N,K-k}) - \sum_{d=1}^J \lambda_d \partial_{x_d} \tilde{f}^{N,K-k} \right) \\ + \left(\frac{\mu}{\sqrt{2}} + \frac{\theta}{\sqrt{\eps}} \right) \sum_{k=0}^K S^{k}_{\delta t} \left(\mathcal{B}(\tilde{f}^{N,K-k}) - \sum_{d=J+1}^{J'} \partial_{x_d} \tilde{f}_d^{N,K-k} \right) + \mathcal{O}((k+1)\eps).
\end{multline}
The difference inside the first brackets in equation \eqref{err_step1} precisely corresponds to the spatial discretization error of the hyperbolic contribution. More specifically, choosing a second order upwind scheme, we can obtain a second order spatial discretization. Unfortunately, if we compute directly the spatial truncation error between the second brackets for the parabolic term using an upwind scheme, a term $\Delta x^2 / \sqrt{\eps}$ appears. Following the idea of \cite{lafitte2012}, an approximation based on centered difference schemes allows us to recover consistency also for the parabolic term. Estimating explicitly $S_{\eps} \Delta_{\mathcal{B}} f$, where $\Delta_{\mathcal{B}} : f^\eps \in C_c^\infty \left( \Omega \times (0,T);\, \R^{L} \right) \longrightarrow \mathcal{B}(f^\eps) -   \sum_{d=J+1}^{J'} \partial_{x_d} \tilde{f}_d^{N,K-k}$. Then, we rewrite
\begin{equation}
    S_{\eps}(\Delta_{\mathcal{B}} f^\eps) = \mathcal{M}\left(\sum_i \Delta_{\mathcal{B}} f_i^\eps \right)-\eps \, \mathcal{B} \,\Delta_{\mathcal{B}}\,f^\eps.
\end{equation}
Consequently, we have 
\begin{equation}
    S^2_{\eps}(\Delta_{\mathcal{B}} f^\eps) = \mathcal{M}\left(\sum_i \Delta_{\mathcal{B}} f_i^\eps \right) - \eps \, \mathcal{M} \left( \sum_i \mathcal{B} \,\Delta_{\mathcal{B}}\,f_i^\eps \right) - \eps \mathcal{B} \mathcal{M}\left(\sum_i \Delta_{\mathcal{B}} f_i^\eps \right) + \eps^2 \mathcal{B}^2 \Delta_{\mathcal{B}} f^\eps,
\end{equation}
which yields, for $k \geq 3$
\begin{multline}
    S^k_{\eps}(\Delta_{\mathcal{B}} f^\eps) = \mathcal{M}\left(\sum_i \Delta_{\mathcal{B}} f_i^\eps \right) - \eps \, \mathcal{M} \left(\sum_i \mathcal{B} \,\Delta_{\mathcal{B}}\,f_i^\eps \right) - \eps \mathcal{B} \mathcal{M}\left(\sum_i \Delta_{\mathcal{B}} f_i^\eps \right) \\
    +\eps^2 \mathcal{M} \left(\sum_i \mathcal{B}^2 \Delta_{\mathcal{B}} f_i^\eps \right) + \eps^2 \mathcal{B}^2 \mathcal{M} \left(\sum_i \Delta_{\mathcal{B}} f_i^\eps \right) \\ + \Phi \mathcal{M} \left( \sum_i \mathcal{B} \Delta_{\mathcal{B}} f_i \right) + \eps^2 (k-3) \mathcal{M}\left(\mathcal{B} \mathcal{M} \left( \sum_i \Delta_{\mathcal{B}} f_i \right) \right) + \mathcal{O}(\eps^3).
\end{multline}
Let us now compute the average value of $S^k_\eps(\Delta_{\mathcal{B}} f^\eps)$, using the fact that $\Delta_{\mathcal{B}}$ is odd.\\
For $k=0$, 
\begin{equation}
    \sum_i \Delta_{\mathcal{B}} f_i^\eps = \eps \sum_i \Delta_{\mathcal{B}} g_i^\eps = \eps \mathcal{O} ( \Delta x^2).
\end{equation}
For $k=1$,
\begin{equation}
    \sum_i S_\eps (\Delta_{\mathcal{B}} f_i^\eps) = \eps \sum_i \Delta_{\mathcal{B}} g_i^\eps - \eps \sum_i \mathcal{B} \Delta_{\mathcal{B}} g_i^\eps = \eps \mathcal{O} ( \Delta x^2) + \eps^2 \mathcal{O}(\Delta x^2).
\end{equation}
More generally, for $k \geq 2$, we obtain
\begin{equation}
    \sum_i S^k_\eps (\Delta_{\mathcal{B}} f_i^\eps) = \eps \mathcal{O} ( \Delta x^2) + \eps^2 \mathcal{O}(\Delta x^2) + \mathcal{O}(\eps^3).
\end{equation}
The last estimate together with Young's inequality leads to
\begin{equation}
     \sum_i e_i^{N, k+1} = \mathcal{O} \left(\sqrt{\eps} (k+1) \Delta x^2 \right) + \mathcal{O} \left((k+1) \eps^{5/2} \right), \qquad \forall \, k \geq 1.
\end{equation}
This implies that
\begin{equation}
    \begin{aligned}
    E^{N+1} = &\frac{\tilde{u}^{N+1} - u^{N,K+1} + \eps \displaystyle \sum_i e^{N, K+1}_i}{\Delta t} - \left(1 - (K+1) \frac{\eps}{\Delta t} \right) \frac{\tilde{u}^{N, K+1} - \tilde{u}^{N,K}}{\eps} \\[7pt]
    - &\left(1 - (K+1) \frac{\eps}{\Delta t} \right) \sum_i \left(e^{N, K+1}_i - e^{N,K}_i \right) \\[7pt]
    = &\left(1 - (K+1) \frac{\eps}{\Delta t} \right) \partial_t u(t^{N, K+1}) + \mathcal{O}(\Delta t) + \eps^{1/2} \mathcal{O} \left( \frac{\Delta x^2 + \eps^2}{\Delta t} \right) \\[7pt]
    - &\left(1 - (K+1) \frac{\eps}{\Delta t} \right) \partial_t u(t^{N,K+1}) + \mathcal{O}(\eps) + \mathcal{O} \left(\frac{\eps^2}{\Delta t} \right) \\
    + &\mathcal{O} (\eps^{1/2} \, \Delta x^2) + \mathcal{O}(\eps^{5/2}) + \mathcal{O} \left( \frac{ \eps^{3/2} \, \Delta x^2}{\Delta t} \right) + \mathcal{O} \left( \frac{ \eps^{7/2}}{\Delta t} \right) 
    \end{aligned}
\end{equation}
In conclusion, we summarize the above results on consistency of the method:
\begin{thm}
    Under the CFL condition $\Delta t = \mathcal{O}(\Delta x^2)$, the PFE scheme with second order upwind flux for $f_i$ with $i=1,\dots,J$ and centered flux for $f_j$ with $j=J+1, \dots, L$ is consistent with \eqref{system_diffusive} at order $1$ in $\eps$ and, as $\eps \to 0$, the limiting scheme is consistent at order $1$ in time and $2$ in space, namely
    \begin{equation}
        E^{N+1} = \mathcal{O}(\eps) + \mathcal{O}(\Delta t)+\mathcal{O}\left(\frac{\eps^2}{\Delta t}\right).
    \end{equation}
\end{thm}

\section{Numerical Simulations}
\label{Section::Numerical_Simulations}

\subsection{Order of convergence: linear diffusion equation}
\begin{figure}[h]
    \begin{minipage}[h!]{0.49\linewidth}
    \begin{center}
    \includegraphics[width=\linewidth]{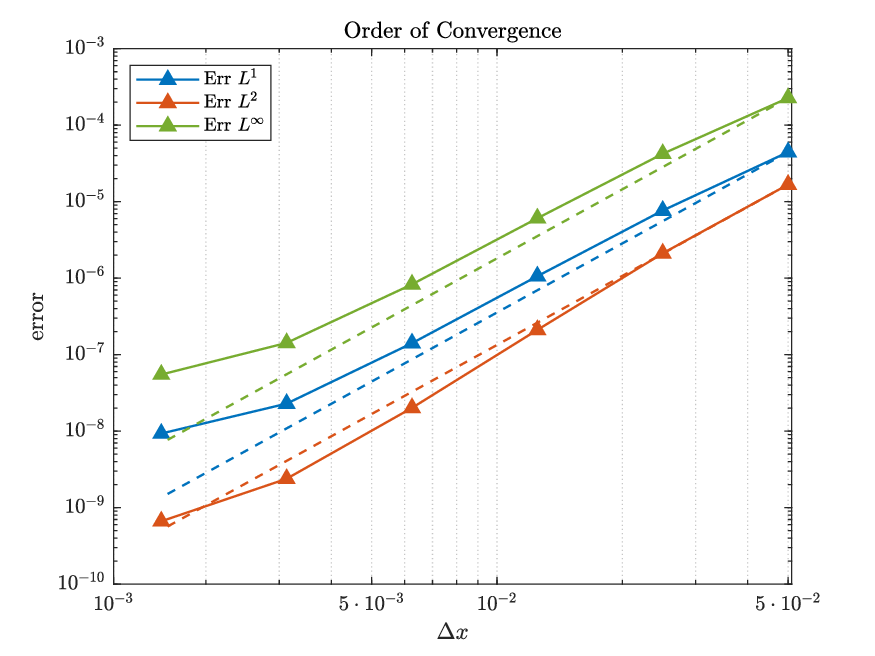}  
    \end{center}
    \end{minipage}
    \hfill
    \begin{minipage}[h]{0.49\linewidth}
    \begin{center}
    \includegraphics[width=\linewidth]{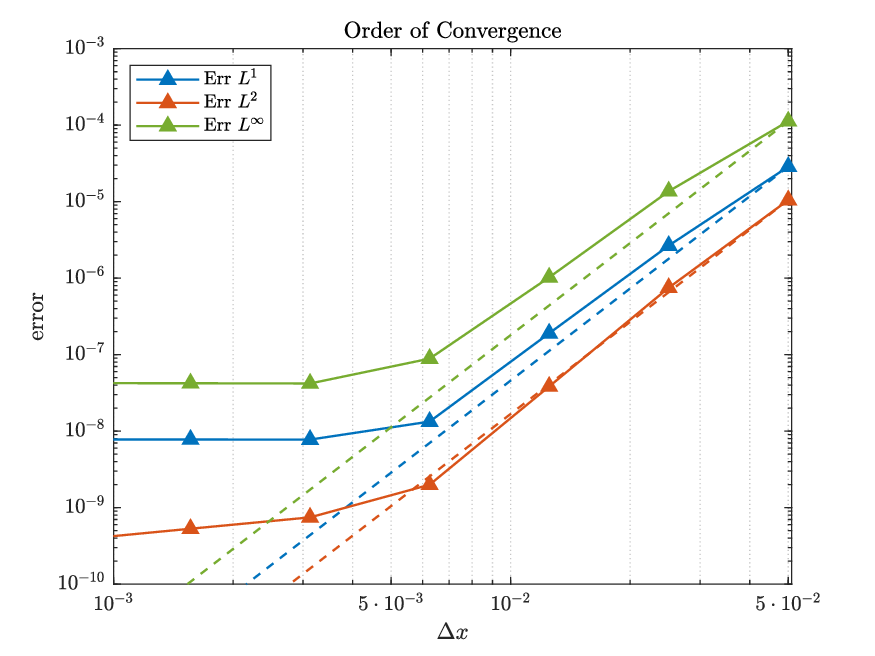} 
    \end{center} 
    \end{minipage}
    \caption{Order of convergence for the $3$rd order scheme (on the left) and the $4$th order scheme (on the right) for different norms, using the linear diffusion equation \eqref{linear_diffusion_test}.}
    \label{fig::spatial_order_linear}
\end{figure}

\begin{figure}[h!]
    \begin{minipage}[h]{0.49\linewidth}
    \begin{center}
    \includegraphics[width=\linewidth]{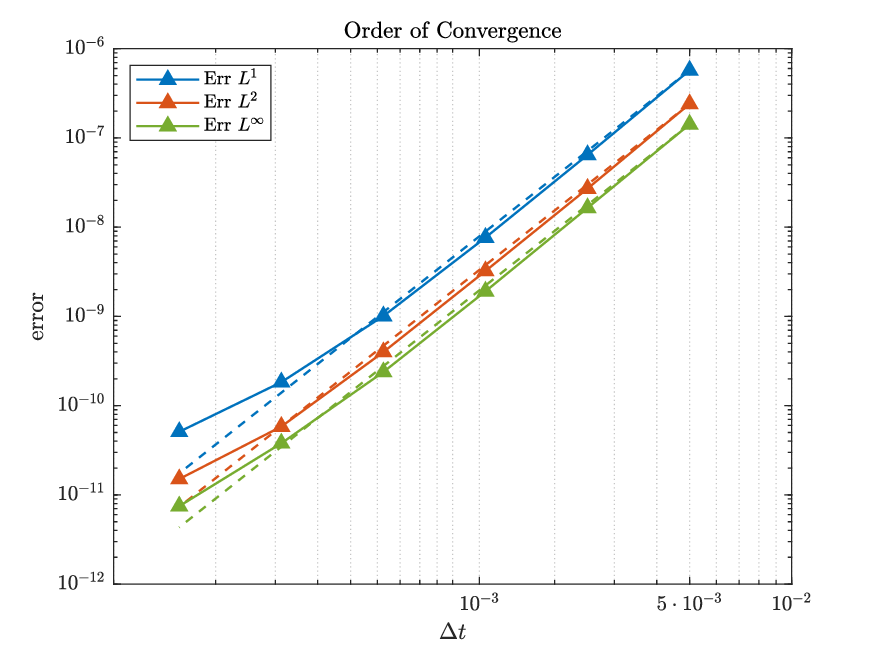}  
    \end{center}
    \end{minipage}
    \hfill
    \begin{minipage}[h]{0.49\linewidth}
    \begin{center}
    \includegraphics[width=\linewidth]{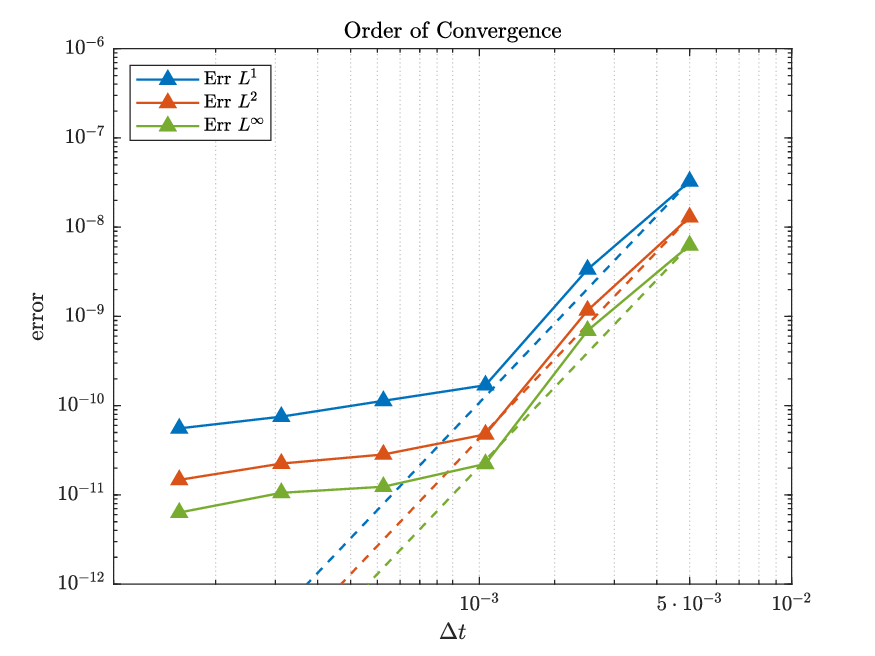} 
    \end{center} 
    \end{minipage}
    \caption{Order of convergence for the PRK3 (on the left) and the PRK4 (on the right) for different norms, using the linear diffusion equation \eqref{linear_diffusion_test}.}
    \label{fig::temporal_order_linear}
\end{figure}

The test for the order of convergence is first performed for the DRM1 kinetic model by considering the linear diffusion problem \cite{wissocqabgrall2024_newlocal} 
\begin{equation}
\label{linear_diffusion_test}
    \partial_t u = \xi \partial_{xx} u,
\end{equation}
where $\xi = 10^{-2}$ and the initial condition is $u_0$, on $(0,T) \times \Omega = (0, T) \times [0,\,1]$. In particular, we consider as initial condition 
\begin{equation}
\label{initial_condition_linear_diffusion}
    u_0(x) = 1 + 0.01 \exp \left(-\frac{(x-0.5)^2}{\delta^2} \right)
\end{equation}
where $\delta=0.1$ and with periodic boundary conditions. 
We compare the numerical solution with the exact solution
\begin{equation}
    u(x,t) = 1 + 0.01 \sqrt{ \frac{\delta^2}{\delta^2 + 4\xi\,t}} \, \exp \left(-\frac{(x-0.5)^2}{\delta^2+4\xi\,t} \right).
\end{equation}
We first test the spatial order of convergence, where the time step is taken equal to $\Delta t = 10^{-7}$ and $\delta t = \eps = 10^{-10}$. With this choice, we are able to study the decrease of the error with respect to the spatial step size $\Delta x$ (since the truncation error dueto the temporal discretization becomes negligible). In Fig \ref{fig::spatial_order_linear} we compare the order of convergence in $L^1$, $L^2$ and $L^\infty$ norm of the $3$rd order scheme and $4$th order scheme, plotting the computed error against the expected error in log-log scale. The truncation error behaves as $\mathcal{O}(\Delta x^3)$ and $\mathcal{O}(\Delta x^4)$, respectively, for large $\Delta x$, until a plateau is reached. Indeed, for small values of $\Delta x$, the contribution of the inner integrator $\mathcal{O}(\eps)$ becomes dominant. This is clearly seen in the $4$th order approximation, where the error does not decrease more than $10^{-10}$, which corresponds to $\eps$ in this test case.\\ 
The temporal order of projective integration method is proven using a different approach. Indeed, since we have to satisfy the parabolic CFL condition $\Delta t = \mathcal{O}(\Delta x^2)$, the contribution of the spatial discretization error is dominant. Our idea is based on Richardson extrapolation \cite{richardson1911, koellermeiersamaey2025}, where two solution computed with different time step $\frac{\Delta t}{k}$, for $k \in \N$ are compared.\\
Assuming a time discretization of order $p$, the solution $u^{n+1}_{\Delta t/k}$ is deviating from the exact solution $u^{n+1}_*$ of the semidiscrete system \eqref{semidiscrete_operator} according to 
\begin{equation}
    u^{n+1}_{\Delta t / k} = u^{n+1}_* + k\, C_n \, \left( \frac{\Delta t}{k} \right)^{p+1} + \mathcal{O} \left( \Delta t^{p+2} \right),
\end{equation}
where $C_n$ is the error constant.
Therefore, we can rewrite the solution $u^{n+1}_*$ as
\begin{equation}
    u^{n+1}_* = \frac{u^{n+1}_{\Delta t} - 2^p u^{n+1}_{\Delta t/2}}{1-2^p},
\end{equation}
from which we obtain an estimation of $u^{n+1}_{\Delta t}-u^{n+1}_*$ in terms of $u^{n+1}_{\Delta t}-u^{n+1}_{\Delta t/2}$, namely
\begin{equation}
    u^{n+1}_{\Delta t}-u^{n+1}_* = \frac{2^p \left( u^{n+1}_{\Delta t}-u^{n+1}_{\Delta t/2} \right)}{2^p-1}.
\end{equation}
Since to compute the order of convergence we are interested in the ratio between two ``consecutive" errors, i.e. errors computed with $\Delta t/k$ and $\Delta t/(2k)$, we can give the following definition of the relative error for a general norm at time $T$ 
\begin{equation}
\label{error_temporal}
    e_{\Delta t} = || u_{\Delta t}(T) - u_{\Delta t/2}(T)||
\end{equation}
One observes that the truncation error behaves as expected and the results are depicted in Fig. \ref{fig::temporal_order_linear}, where plot the slope of the expected error n log-log scale. In this case, as for the spatial discretization error, the truncation error reaches a plateau, since the contribution of the inner integrator $\mathcal{O}(\eps)$ becomes dominant.

\subsection{Order of convergence: advection-diffusion equation}
\begin{figure}[h!]
    \begin{minipage}[h]{0.49\linewidth}
    \begin{center}
    \includegraphics[width=\linewidth]{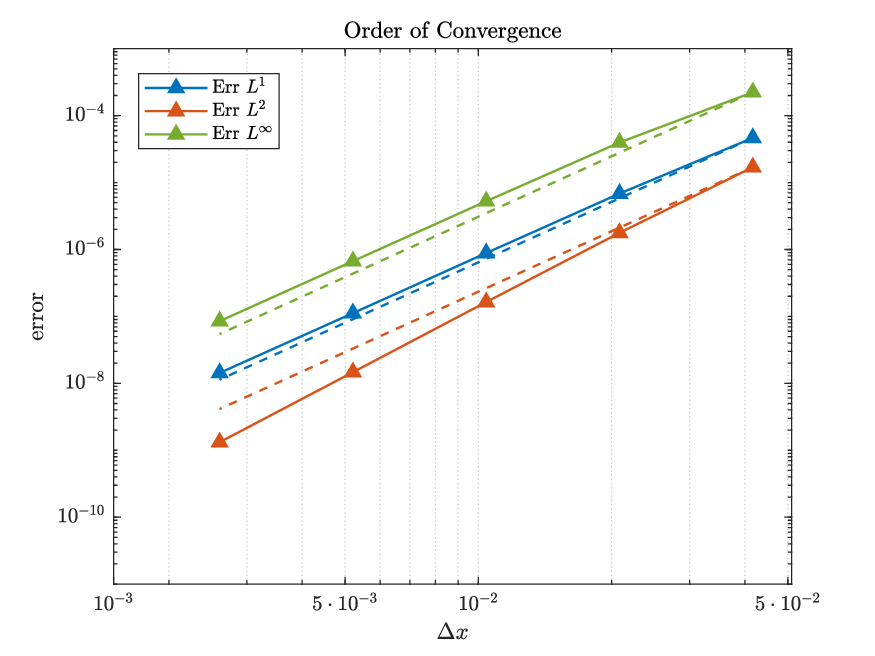}  
    \end{center}
    \end{minipage}
    \hfill
    \begin{minipage}[h]{0.49\linewidth}
    \begin{center}
    \includegraphics[width=\linewidth]{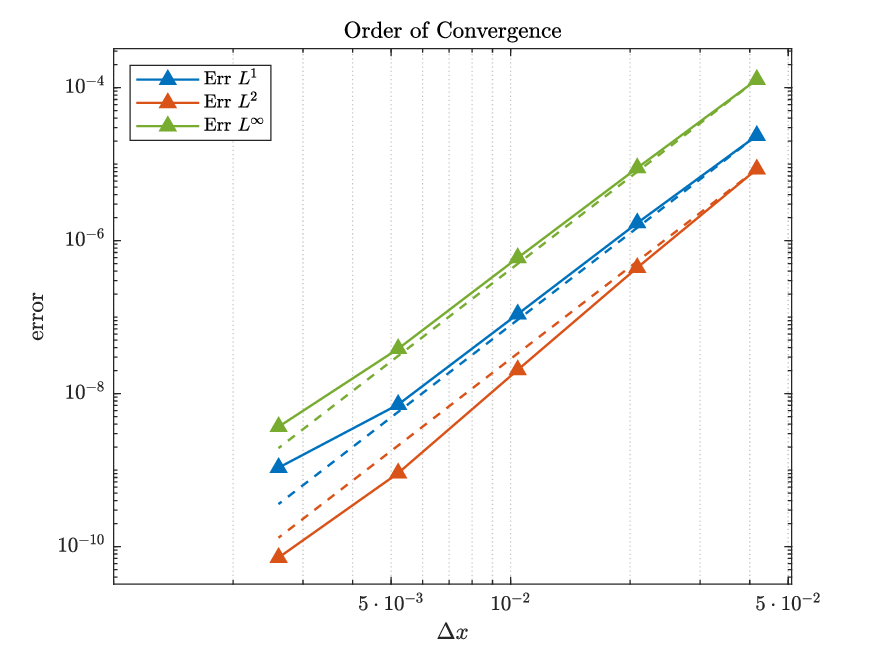} 
    \end{center} 
    \end{minipage}
    \caption{Order of convergence for the $3$rd order scheme (on the left) and the $4$th order scheme (on the right) for different norms, using the advection-diffusion equation \eqref{linear_diffusion_test}.}
    \label{fig::spatial_order_test2}
\end{figure}

Following again \cite{wissocqabgrall2024_newlocal}, we propose a test for the advection-diffusion equation with constant velocity $c$
\begin{equation}
    \partial_t u + c \partial_x u = \partial_{xx} u,
\end{equation}
with the same initial condition as \eqref{initial_condition_linear_diffusion}. The choice for the discretization parameters is analogous to the previous test. Setting $c=10$, we reproduce the same test case, where now the solution is advected with a positive velocity. 
We first test again the spatial order of convergence, in the same setting of the previous case. In Fig. \ref{fig::spatial_order_test2} we compare again the order of convergence in $L^1$, $L^2$ and $L^\infty$ norm of the $3$rd order and $4$th order schemes, plotting the computed error against the expected error in log-log scale.

\begin{figure}[h!]
    \begin{minipage}[h]{0.49\linewidth}
    \begin{center}
    \includegraphics[width=\linewidth]{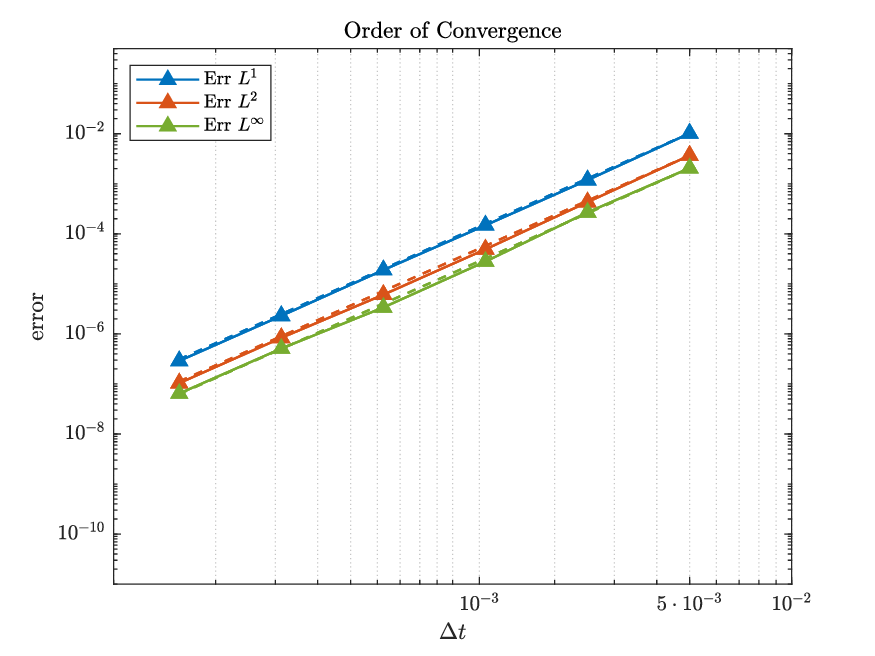}  
    \end{center}
    \end{minipage}
    \hfill
    \begin{minipage}[h]{0.49\linewidth}
    \begin{center}
    \includegraphics[width=\linewidth]{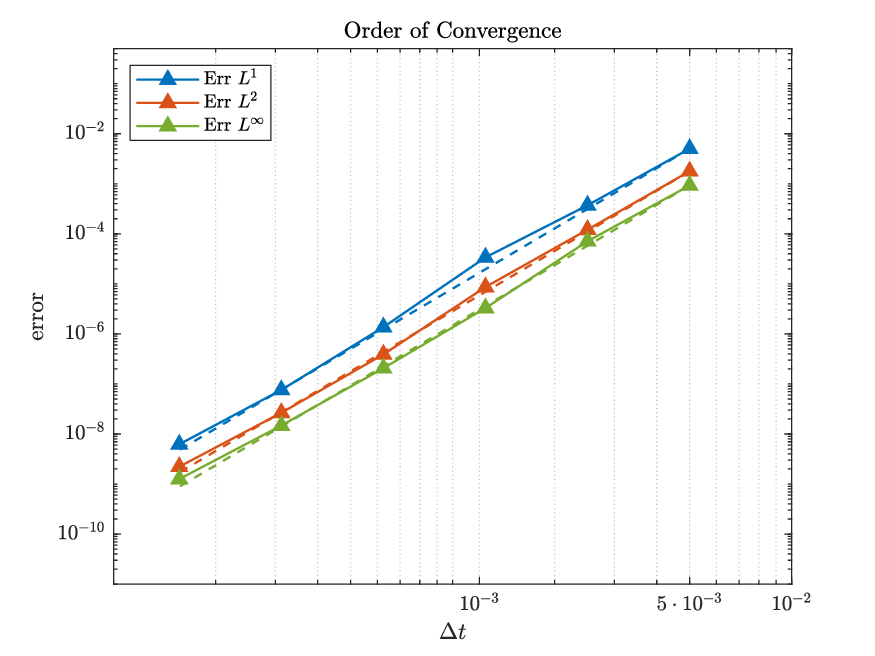} 
    \end{center} 
    \end{minipage}
    \caption{Order of convergence for the PRK3 (on the left) and the PRK4 (on the right) for different norms, using the advection-diffusion equation \eqref{linear_diffusion_test}.}
    \label{fig::temporal_order_test2}
\end{figure}

The temporal order of convergence is tested again using the definition in \eqref{error_temporal}. In Fig. \ref{fig::temporal_order_test2} we report the slope of the expected error in log-log scale. Differently from the previous test, the plateau is not reached here, since the spatial or the temporal error is dominant compared to the contribution of the inner integrator $\mathcal{O}(\eps)$.

\subsection{Testing the speedup: the viscous LWR equation}
\begin{figure}[h!]
    \begin{minipage}[h]{0.49\linewidth}
    \begin{center}
    \includegraphics[width=\linewidth]{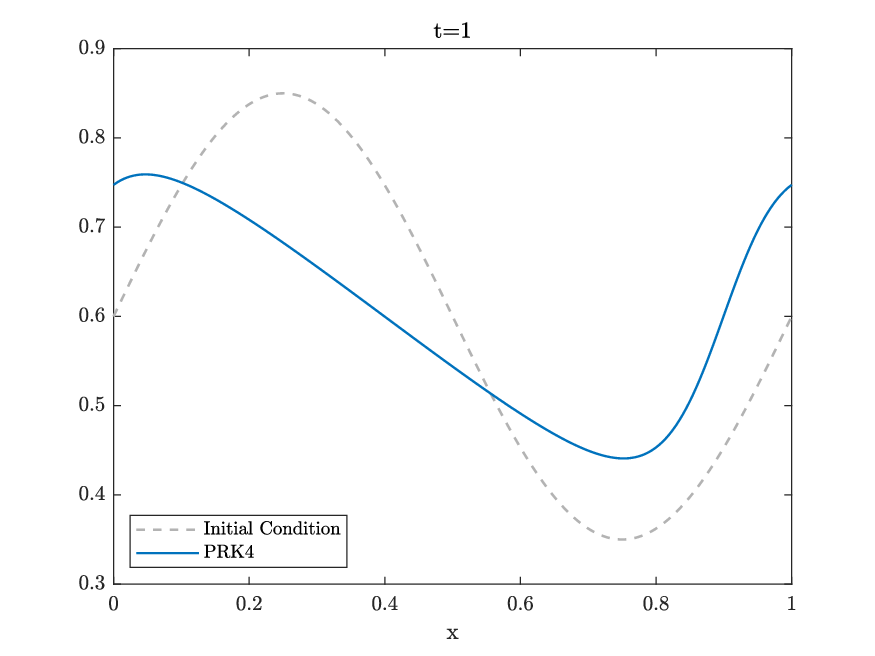}  
    \end{center}
    \end{minipage}
    \hfill
    \begin{minipage}[h]{0.49\linewidth}
    \begin{center}
    \includegraphics[width=\linewidth]{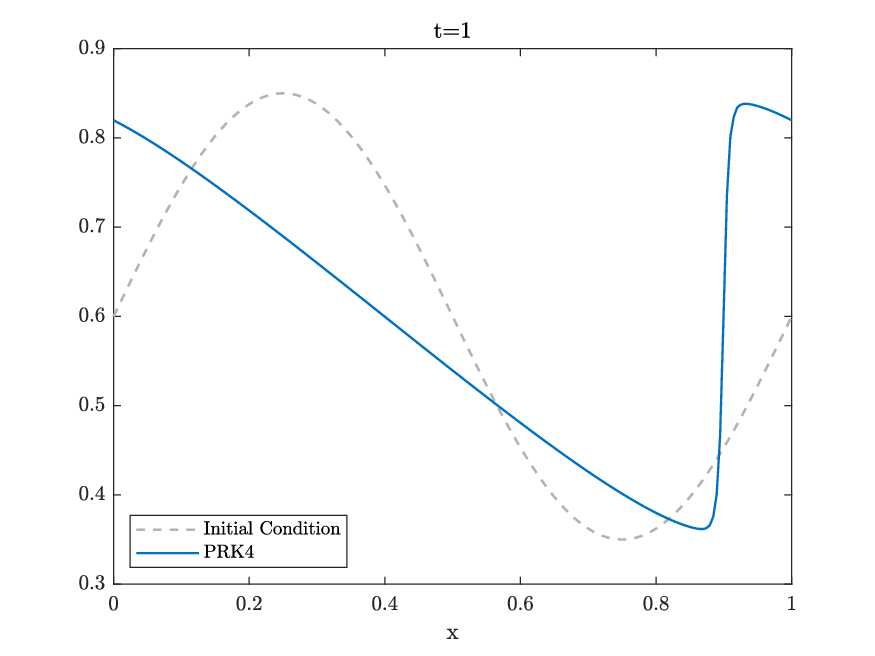} 
    \end{center} 
    \end{minipage}
    \caption{Numerical approximation of the viscous LWR equation for two different values of $\xi$, respectively $\xi=10^{-2}$ (on the left) and $\xi=10^{-3}$ (on the right).}
    \label{Test_ViscousLWR}
\end{figure}

Let us consider the viscous LWR equation
\begin{equation}
    \partial_t u + \partial_x \left(\frac{u-u^2}{2} \right) = \xi \partial_{xx}\,u,
\end{equation}
where $\xi$ assumed different values. We consider the following initial condition, as in \cite{ferretti2025}
\begin{equation}
    u_0(x) = 0.6 + 0.25\,\sin(2\,\pi\,x),
\end{equation}
in the domain $[0,1]$ with periodic boundary conditions. The discretization is performed using $\Delta x = 5 \cdot 10^{-3}$ and $\Delta t$ according to the parabolic CFL condition. The projective integration parameters are $\delta t= \eps$ and $K=2$, where $\eps=10^{-7}$ from now on. In Fig. \ref{Test_ViscousLWR} we report the solution for $\xi=10^{-2}$ and $\xi = 10^{-3}$. Let us now test the speedup of the numerical approximation for this test case. From a theoretical point of view, if we assume that the overhead due to extrapolation si negligible and time stepping with the innermost integrator is computationally most demanding, the speedup can be computed as 
\begin{equation}
    \mathcal{S}_{\text{PI}} = \frac{\Delta t}{\delta t (K+1)},
\end{equation}
that is the ratio of the total number of naive forward Euler time steps within one outermost time step over the number of real innermost steps in the method. In Table \ref{Table::CPU_improvement} we report the CPU time of the projective integration method, compared to the classical explicit integration. When $\eps \to 0$, the projective integration approach becomes much more efficient than the direct integration. Moreover, in Table \ref{Table::Comparison_IMEX} we compare a PFE method with a first order IMEX scheme for stiff BGK equation proposed in \cite{pareschirusso2005}. The numerical scheme based on projective integration method exhibits a computational cost of the same order of the IMEX approach (although slightly slower) and it proves to be \textit{asymptotic preserving} with respect to $\eps$. The great advantage of the projective integration approach is that high order extension though a PRK scheme \eqref{PRK_scheme}-\eqref{PRK_scheme_2} is easily constructed.

\begin{table}[h!]
    \centering
    \begin{tabular}{c|c|c|c|c|}
         $\eps$ & \multicolumn{2}{c|}{CPU time (s)} & \multicolumn{2}{c|}{Improvement factor} \\
         \hline
         & Direct & Projective Integration &  Real & Theoretical \\
         \hline
         
         $10^{-5}$ & $1.51$ & $0.89$ & $1.69$ & $1.70$ \\
         \hline
         $10^{-7}$ & $152$ & $0.94$ & $162$ & $171$ \\
    \end{tabular}
    \caption{CPU time and improvement factor}
    \label{Table::CPU_improvement}
\end{table}

\begin{table}[h!]
    \centering
    \begin{tabular}{c|c|c|}
         $\eps$ & \multicolumn{2}{c|}{CPU time (s)} \\
         \hline
        & Projective Integration &  IMEX Scheme \\
        \hline
         $10^{-5}$ & $0.89$ & $0.73$ \\
         \hline
         $10^{-7}$ & $0.94$ & $0.62$ \\
    \end{tabular}
    \caption{Comparison between the performance of PFE and IMEX schemes.}
    \label{Table::Comparison_IMEX}
\end{table}

\subsection{The viscous Burgers' equation}
\begin{figure}
    \centering
    \includegraphics[width=0.75\linewidth]{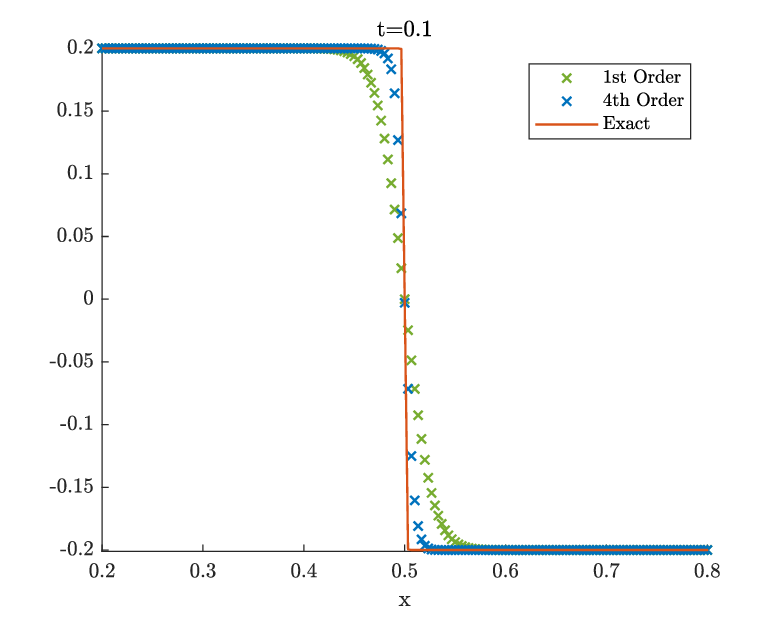}
    \caption{Comparison of $1st$ order and $4th$ order approximation of the ``steady shock'' solution to the viscous Burgers' equation \eqref{Steady_Viscous}.}
    \label{fig:visc_burg1}
\end{figure}

Let us consider the viscous Burgers' equation
\begin{equation}
    \partial_t u + \partial_x \left(\frac{1}{2} \, u^2 \right) = \xi \partial_{xx}\,u,
\end{equation}
where $xi=10^{-3}$. We first test our method for the ``steady shock'' whose exact solution is given by \cite{benton1972, wissocqabgrall2024}
\begin{equation}
\label{Steady_Viscous}
    u(x) = - \frac{2\,\xi}{\delta}\,\tanh \left( \left(x-0.5 \right)/\delta \right),
\end{equation}
where $\delta= 0.01$. In Fig. \ref{fig:visc_burg1} we show the results obtained as numerical approximation with schemes of different order. The discretization is performed with step size $\Delta x = 3 \cdot 10^{-3}$ and $\Delta t$ is chosen according to the parabolic CFL condition. The projective integration parameters are chosen as $\delta t= \eps$ and $K=2$ projective levels. Compared to first-order approximation, fourth order PRK4 scheme is able to reduce numerical dissipation errors, preserving the steady state. The simulation has a perfect agreement with \cite{wissocqabgrall2024}.

\subsection{The viscous Burgers' equation with strongly degenerate diffusion}
\begin{figure}[h!]
    \begin{center}
    \includegraphics[width=0.85\linewidth]{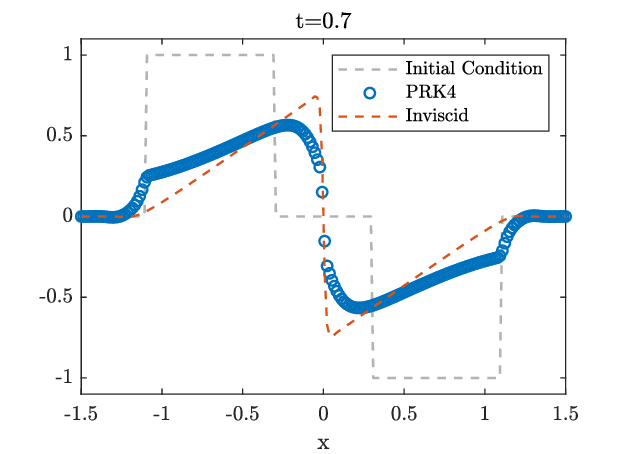}  
    \end{center}
    \caption{Numerical approximation of the viscous Burgers' equation compared with the inviscid solution at final time $t=0.7$.}
    \label{Test3_ViscousBurgers1D}
\end{figure}

Let us consider the Burgers equation 
\begin{equation}
\label{Burgers_Diffusion_1D}
    \partial_t u + \partial_x (u^2) = \xi \partial_x\left(\nu(u) \partial_x u\right),
\end{equation}
with strongly degenerate diffusion
\begin{equation}
\label{eq:strongly_degenerate_burgers_1D}
    \nu(u) = \begin{cases}
        0, \quad &|u| \leq 0.25,\\
        1, \quad &|u| > 0.25,
    \end{cases}
\end{equation}
and the following initial condition
\begin{equation}
    u_0(x) = \begin{cases}
        1, \qquad &\text{if } \displaystyle \displaystyle -\frac{1}{\sqrt{2}} - 0.4 \leq x \leq \displaystyle -\frac{1}{\sqrt{2}} + 0.4,\\[7pt]
        -1, \qquad &\text{if } \displaystyle \frac{1}{\sqrt{2}} - 0.4 \leq x \leq \displaystyle \frac{1}{\sqrt{2}} + 0,\\[7pt]
        0, \qquad &\text{otherwise}.
    \end{cases}
\end{equation}
We set $\xi=0.1$ and the discretization parameters are chosen $\Delta x = 0.015$, $\delta t=\eps$, $\Delta t = C\,\Delta x^2$ and $K=2$. In Fig. \ref{Test3_ViscousBurgers1D} we show the numerical results at final time $t=0.7$, compared with the inviscid case ($\xi=0$). This test shows the ability of the scheme to deal with strong degenerate diffusions.

\subsection{The two-dimensional viscous Burgers' equation}
\begin{figure}[h!]
    \begin{minipage}[h]{0.49\linewidth}
    \begin{center}
    \includegraphics[width=0.95\linewidth]{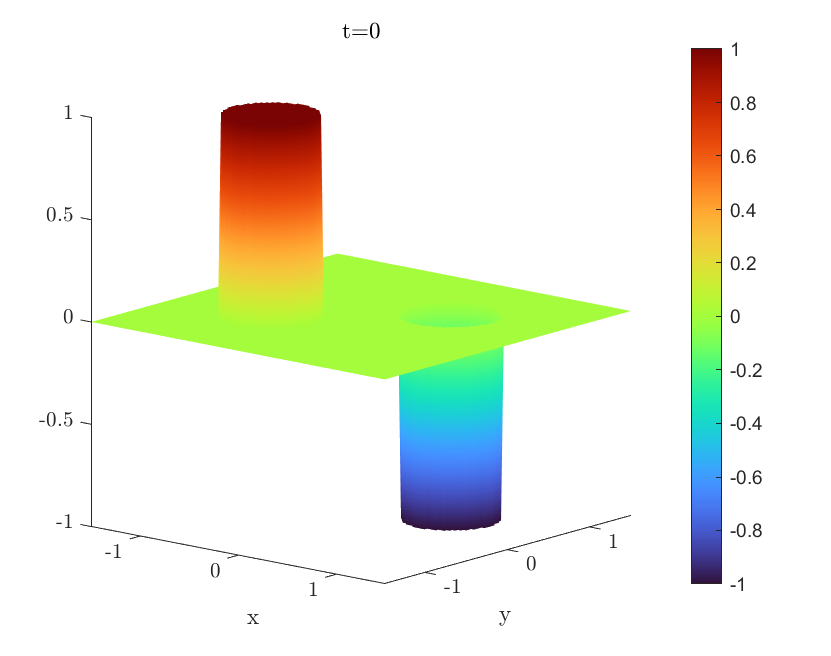}  
    \end{center}
    \end{minipage}
    \hfill
    \begin{minipage}[h]{0.49\linewidth}
    \begin{center}
    \includegraphics[width=0.95\linewidth]{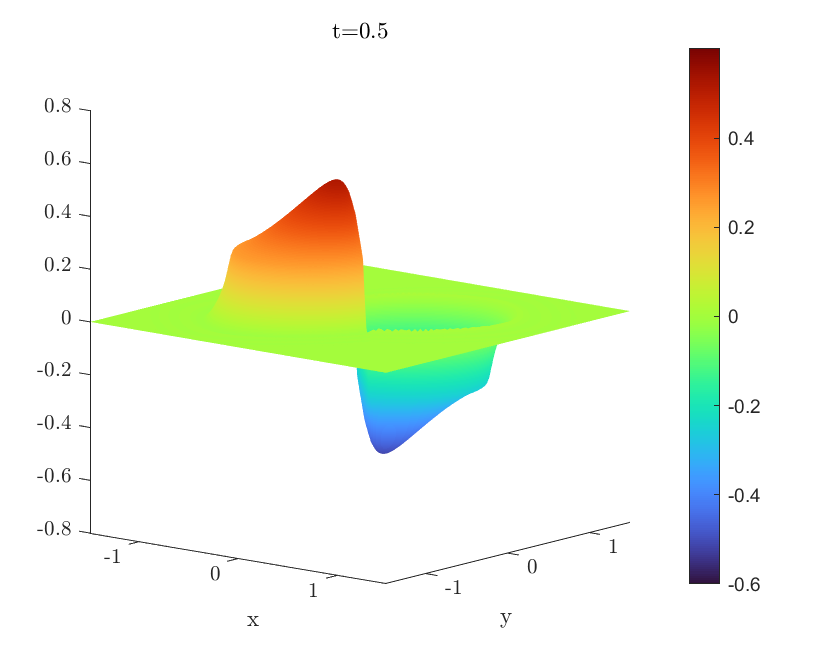} 
    \end{center} 
    \end{minipage}
    \caption{Numerical approximation of the two-dimensional viscous Burgers' equation: initial condition (left) and final solution at time $t=0.5$ (right).}
    \label{Test_2D_Burgers}
\end{figure}

We propose here an extension of \eqref{Burgers_Diffusion_1D} to higher dimensions. Let us focus on the $2$D viscous Burgers' equation, explicitly written as
\begin{equation}
    \partial_t u + \partial_x \left(u^2 \right) + \partial_y \left(u^2 \right) = \xi \partial_x\left(\nu(u) \partial_x u\right) + \xi \partial_y\left(\nu(u) \partial_y u\right),
\end{equation}
again with strongly degenerate diffusion \eqref{eq:strongly_degenerate_burgers_1D} and initial condition
\begin{equation}
    u_0(x,y)= \begin{cases}
        -1 \qquad &\text{if } (x-0.5)^2+(y-0.5)^2 \leq 0.16,\\
        1 \qquad &\text{if } (x+0.5)^2+(y+0.5)^2 \leq 0.16,\\
        0 \qquad &\text{otherwise},\\
    \end{cases}
\end{equation}
The projective integration parameters are $K=2$, $\delta t = \eps$ and $\Delta t = $. The spatial discretization is obtained using $\Delta x = \Delta y = 0.03$ in the domain $[-1.5, 1.5] \times [-1.5, 1.5]$. The simulation in Fig. \ref{Test_2D_Burgers} shows the good approximation with PRK4 scheme. We underline again that a great advantage of discrete kinetic BGK approximation is its feasibility of easily extending the strategy to the multidimensional case without increasing the complexity of the structure of system \eqref{kinetic_approximation}.

\subsection{A three-phase flow model}
\begin{figure}[h!]
    \begin{minipage}[h]{0.9\linewidth}
    \begin{center}
    \includegraphics[width=0.9\linewidth]{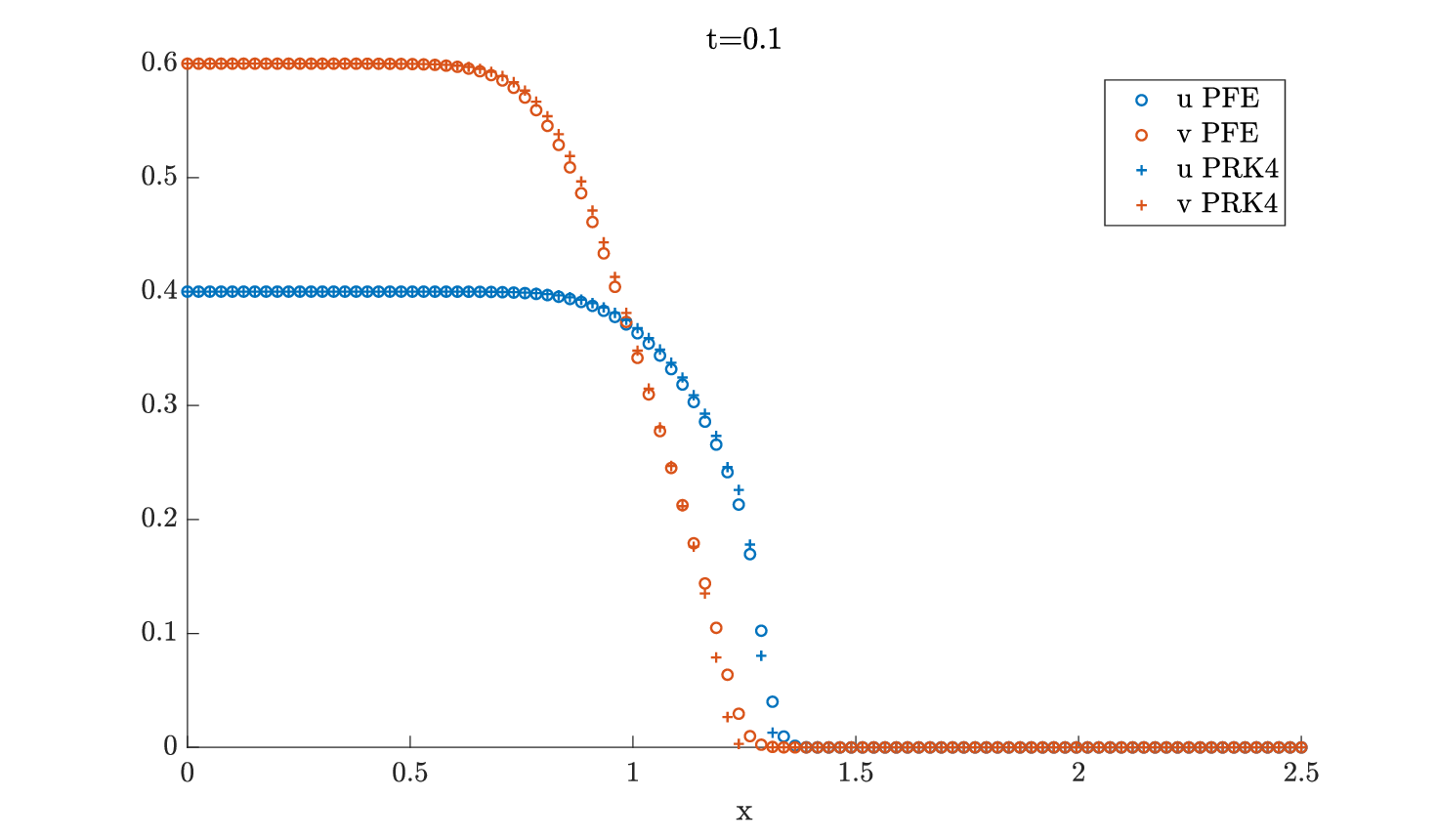}  
    \end{center}
    \end{minipage}
    \vfill
    \begin{minipage}[h]{0.9\linewidth}
    \begin{center}
    \includegraphics[width=0.9\linewidth]{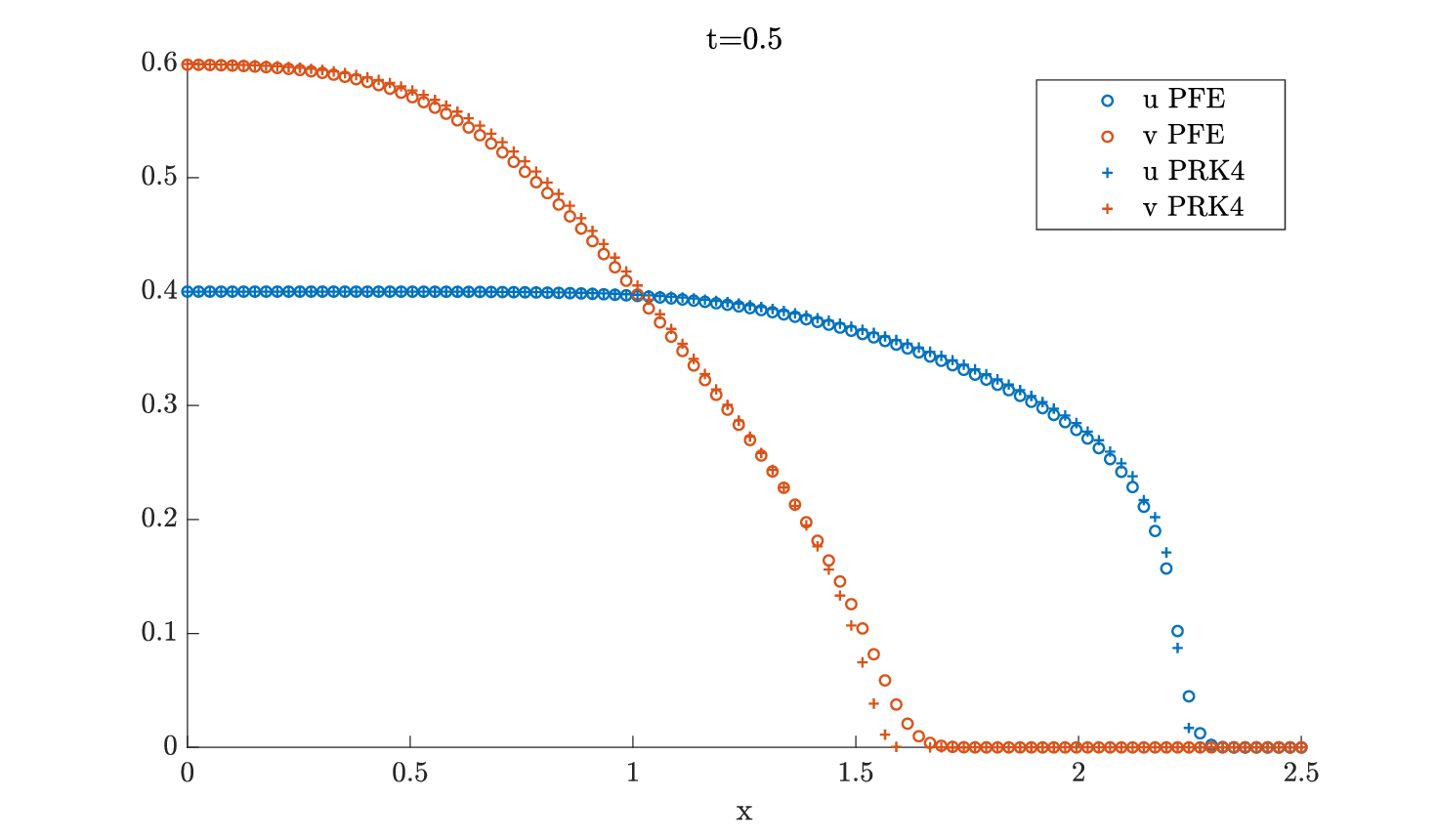} 
    \end{center} 
    \end{minipage}
    \caption{Numerical approximation of the Buckley-Leverett equation with diffusion for the two phase saturations: gas (in blue) and water (in orange).}
    \label{Test1}
\end{figure}

\begin{figure}[h!]
    \begin{minipage}[h]{0.48\linewidth}
    \begin{center}
    \includegraphics[width=0.95\linewidth]{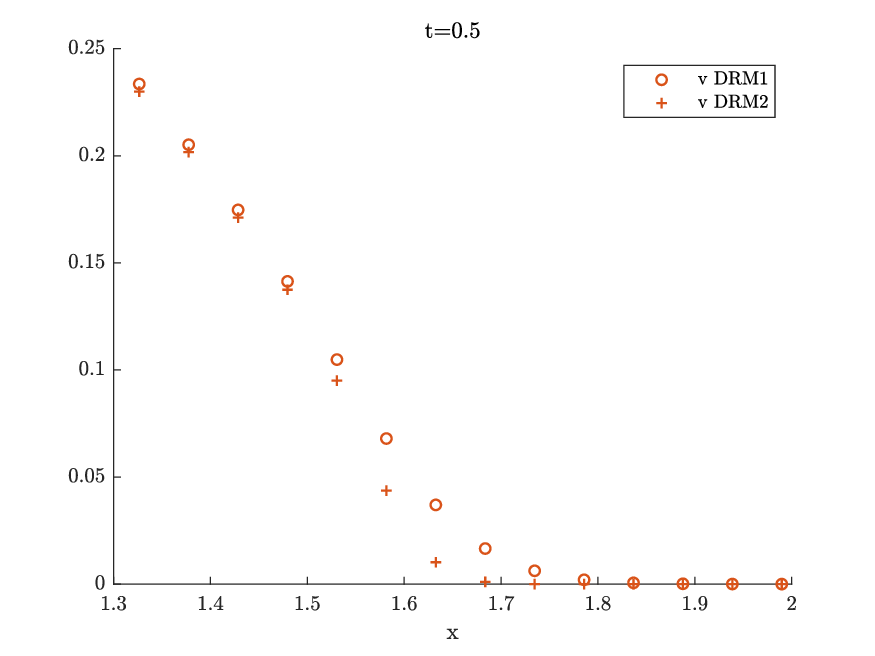}  
    \end{center}
    \end{minipage}
    \hfill
    \begin{minipage}[h]{0.48\linewidth}
    \begin{center}
    \includegraphics[width=0.95\linewidth]{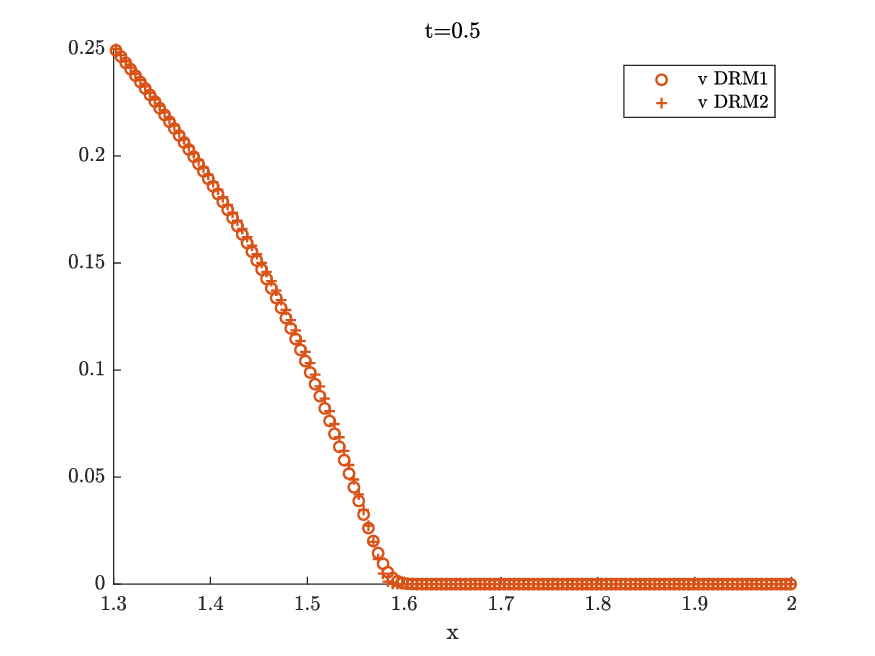} 
    \end{center} 
    \end{minipage}
    \caption{Comparison between the DRM1 (circles) and the DRM2 (crosses) model close to the shock for the Buckley Leverett equation \eqref{BL_Test1_3phase}.}
    \label{Comparison_DRM1DRM2}
\end{figure}

Let us consider the following three phase flow model introduced in the paper of Karlsen et al. \cite{karlsen2001} and approximated also in \cite{aregba2004} to validate the relaxation scheme. It consists of one dimensional system of advection-diffusion equations modeling a multi-phase flow in a porous medium, which reads as follows
\begin{equation}
\label{BL_Test1_3phase}
    \begin{cases}
        \partial_t u + \partial_x f(u) = \xi \partial_{xx} B(u),\\
        \partial_t v + \partial_x g(u, v) = \xi \partial_{xx} B(v),
    \end{cases}
\end{equation}
with initial condition given by 
\begin{equation}
    \left(u_0(x), v_0(x)\right) = \begin{cases}
        (0.4, 0.6), \qquad &\text{if } x < 1,\\
        (0, 0), \qquad &\text{otherwise.}
    \end{cases}
\end{equation}
The variables $u$ and $v$ represent the phase saturations (gas and water, respectively), which evolves according to different dynamics. Indeed, the decoupling between the gas phase and the other phase is a reasonable assumption \cite{karlsen2001} and the flux function $\left(f(u), g(u,v)\right)$ is given by
\begin{equation}
    \begin{cases}
        f(u)&=\displaystyle \frac{u^2}{u^2+(1-u)^2/10},\\[9pt]
        g(u, v)&=\displaystyle \frac{(1-u)^2+u^2/10}{10u^2+(1-u)^2}\,f(v),
    \end{cases}
\end{equation}
and the diffusion function is chosen such that
\begin{equation}
    B'(w)=4w\,(1-w).
\end{equation}
The choice of this particular diffusion function is due to the possibility of reproducing the shape of reservoir models. Moreover, we set the scaling factor $\xi=0.1$ and the relaxation parameter $\eps=10^{-7}$. In Fig. \ref{Test1} we show the numerical results for $\Delta x = 0.005$ in $[0, 2.5]$. The projective integration parameters are chosen as $\delta t= \eps$ and $\Delta t = C\Delta x^2$, with $K=2$ projective levels. The results are in perfect agreement with the one shown in \cite{aregba2004}. We observe that the choice of a $4th$ order scheme leads to a less diffusive approximation of the solution, as expected.\\
We also perform a comparison between DRM1 and DRM2 models, as shown in Fig. \ref{Comparison_DRM1DRM2} with a zoom on the shock. For $\Delta x = 0.05$ the approximation given by the DRM2 model is less dissipative, since the choice of $\lambda_1 \neq -\lambda_2$ improves the kinetic approximation of the original system. Nevertheless, using $\Delta x = 0.0125$, we observe that the two approximations gives equivalent results. 

\subsection{A three-phase flow model with gravitational effect}
\begin{figure}[h!]
    \begin{minipage}[h]{0.49\linewidth}
    \begin{center}
    \includegraphics[width=0.75\linewidth]{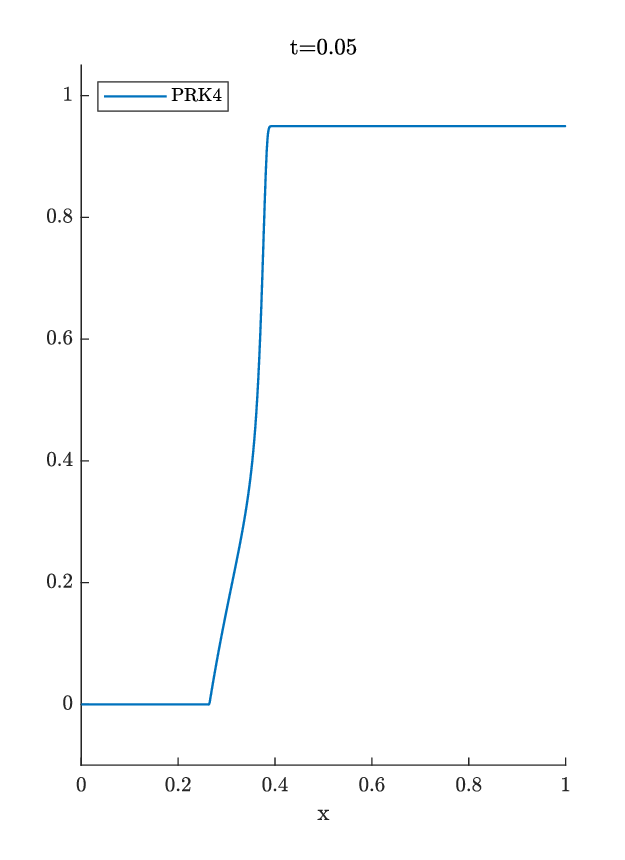}  
    \end{center}
    \end{minipage}
    \hfill
    \begin{minipage}[h]{0.49\linewidth}
    \begin{center}
    \includegraphics[width=0.75\linewidth]{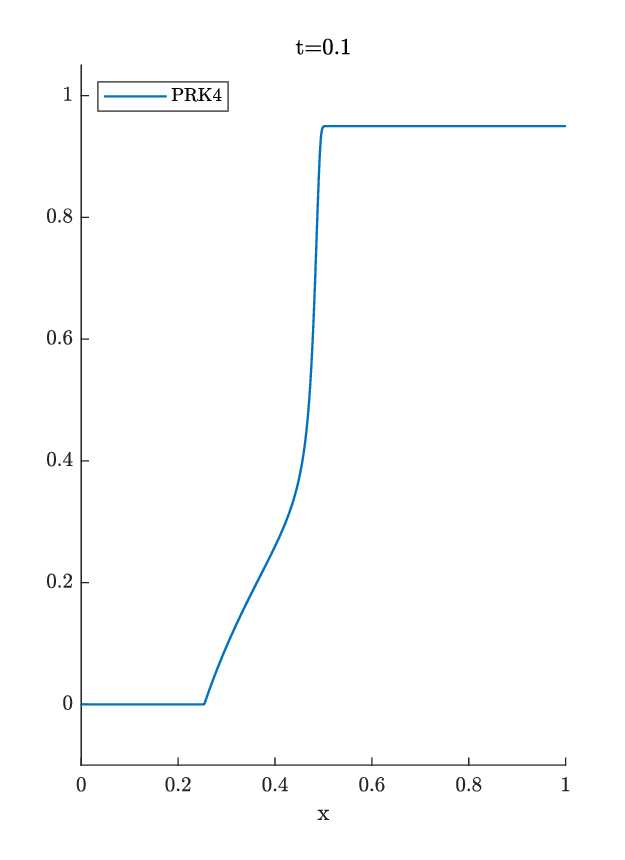} 
    \end{center} 
    \end{minipage}
    \caption{Numerical approximation of the Buckley-Leverett with diffusion \eqref{BL_grav} with $g=0$ at two different times $t=0.05$ and $t=0.1$. }
    \label{Test2_g0}
\end{figure}

\begin{figure}[h!]
    \begin{minipage}[h]{0.49\linewidth}
    \begin{center}
    \includegraphics[width=0.75\linewidth]{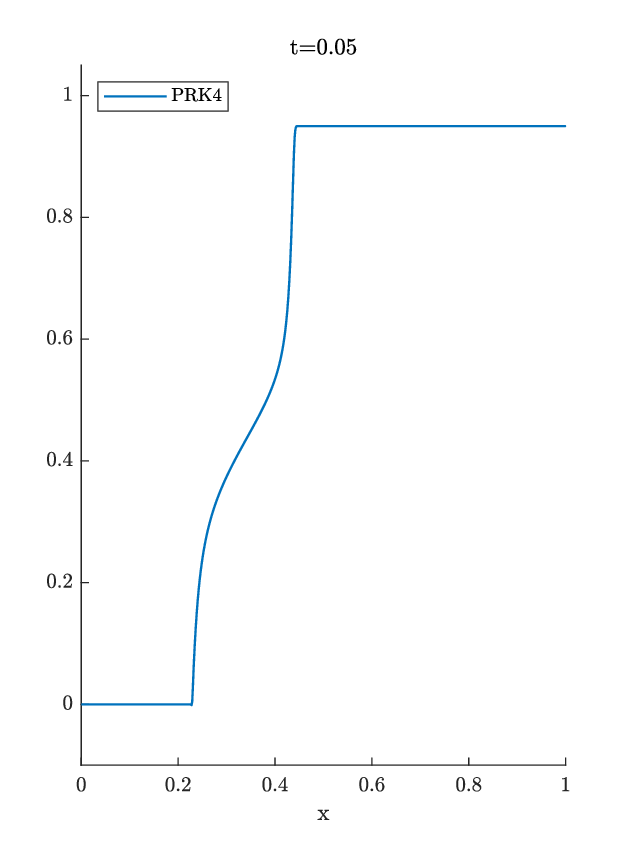}  
    \end{center}
    \end{minipage}
    \hfill
    \begin{minipage}[h]{0.49\linewidth}
    \begin{center}
    \includegraphics[width=0.75\linewidth]{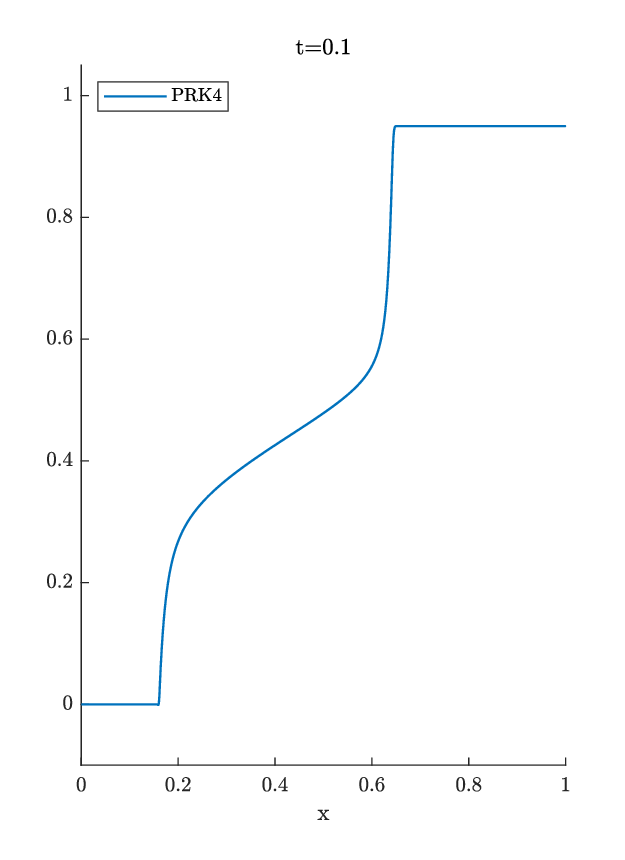} 
    \end{center} 
    \end{minipage}
    \caption{Numerical approximation of the Buckley-Leverett with diffusion \eqref{BL_grav} with $g=5$ at two different times $t=0.05$ and $t=0.1$.}
    \label{Test2_g5}
\end{figure}

We consider again the advection-diffusion Buckley-Leverett equation of the same form 
\begin{equation}
    \partial_t u + \partial_x f(u) = \xi \partial_{xx} B(u),
\end{equation}
where now the flux function is given by
\begin{equation}
    f(u) = \displaystyle \frac{u^2}{u^2 + a\,(1-u)^2} ( 1-g\,(1-u)^2),
\end{equation}
where $a$ is the ratio of the viscosities of the two fluids and $g$ is the gravitational effect. In this test the parabolic term has the form
\begin{equation}
\label{BL_grav}
    B(u) = \begin{cases}
        0, \qquad &\text{if } u <0,\\[7pt]
        \displaystyle \xi \left(2u^2 - \frac{4}{3} u^3 \right) \qquad &\text{if } 0 \leq u \leq 1,\\[7pt]
        \displaystyle \frac{2}{3} \xi \qquad &\text{if } u >1.
    \end{cases}
\end{equation}
Following \cite{kossaczka2022neural}, we set the following initial condition
\begin{equation}
    u_0(x) = \begin{cases}
        0, \qquad &\text{if } 0 \leq x \displaystyle \leq 1-\frac{1}{\sqrt{2}},\\[7pt]
        1, \qquad &\text{if } \displaystyle \frac{1}{\sqrt{2}} < x \leq 1.
    \end{cases}
\end{equation}
We propose some numerical simulations for different values of $g$. We set the scaling factor $\xi=0.01$ and the relaxation parameter $\eps=10^{-7}$. The spatial numerical grid is chosen $\Delta x=0.01$ in $[0,1]$ and the projective integration parameters are $\delta t = \eps$, $\Delta t = C \Delta x^2$ and $K=2$. In Fig. \ref{Test2_g0}  and Fig. \ref{Test2_g5} we observe the results for $g=0$ and $g=5$ respectively at time $t=0.05$ and $t=0.1$.
\subsection{A two-dimensional three-phase model}
As last example, we consider the two-dimensional Buckley-Leverett equation with diffusion, of the form
\begin{equation}
\label{eq_BL_2D}
    \partial_t u + \partial_x f_1(u) + \partial_y f_2(u) = \xi \left(\partial^2_{xx} u + \partial^2_{yy} u \right),
\end{equation}
where the scaling factor $\xi=0.01$ and the flux functions are
\begin{equation}
    f_1(u) = \frac{u^2}{u^2+(1-u)^2}, \qquad f_2(u) = f_1(u)\,\left(1-5(1-u)^2 \right).
\end{equation}
The initial condition is given by
\begin{equation}
    u(x,y,0) = \begin{cases}
        1, \qquad &\text{if } \, x^2+y^2 < 0.5,\\
        0 \qquad &\text{otherwise.}
    \end{cases}
\end{equation}
The domain is $[-1.5, 1.5] \times [-1.5, 1.5]$ and the spatial discretization is obtained taking $\Delta x = \Delta y = 0.015$. The parameters for the projective integration method are the same as before, namely $\delta t = \eps$, $\Delta t = C \Delta x^2$ and $K=2$, where $C=0.4$. Since this simulation has a discontinuous initial condition in two dimensions which may lead to oscillatory behavior, we employ a CWENO3 approximation for the hyperbolic variables, following the scheme \cite{levypupporusso2000}.  The results at time $t=0.5$ agree with the ones of \cite{tadmor2000}. We underline again that the advantage of using relaxation schemes relies on their easily extensibility to the multi-dimensional case, without increasing the complexity of the numerical scheme.

\begin{figure}[h!]
    \begin{minipage}[h]{0.49\linewidth}
    \begin{center}
    \includegraphics[width=0.95\linewidth]{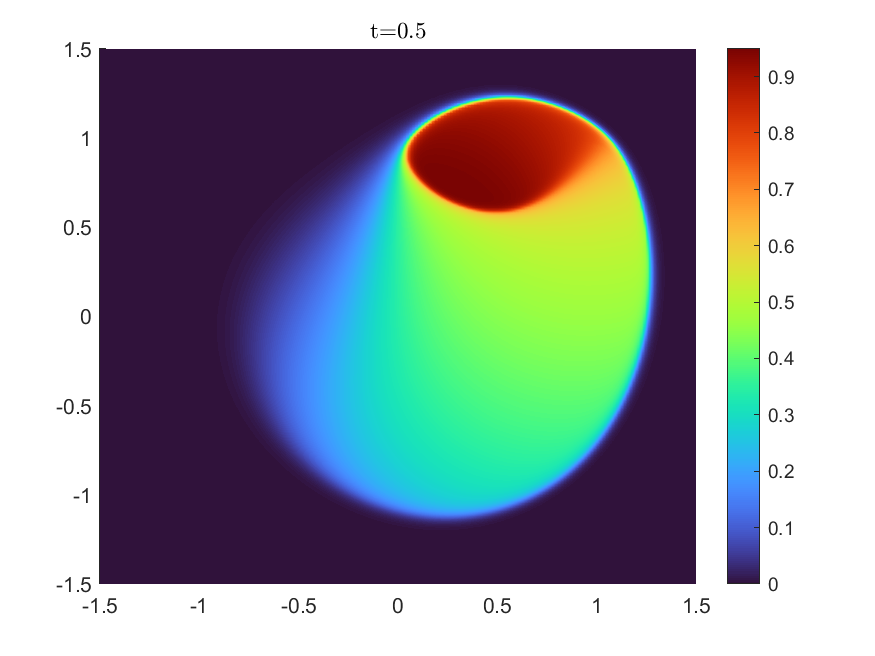}  
    \end{center}
    \end{minipage}
    \hfill
    \begin{minipage}[h]{0.49\linewidth}
    \begin{center}
    \includegraphics[width=0.95\linewidth]{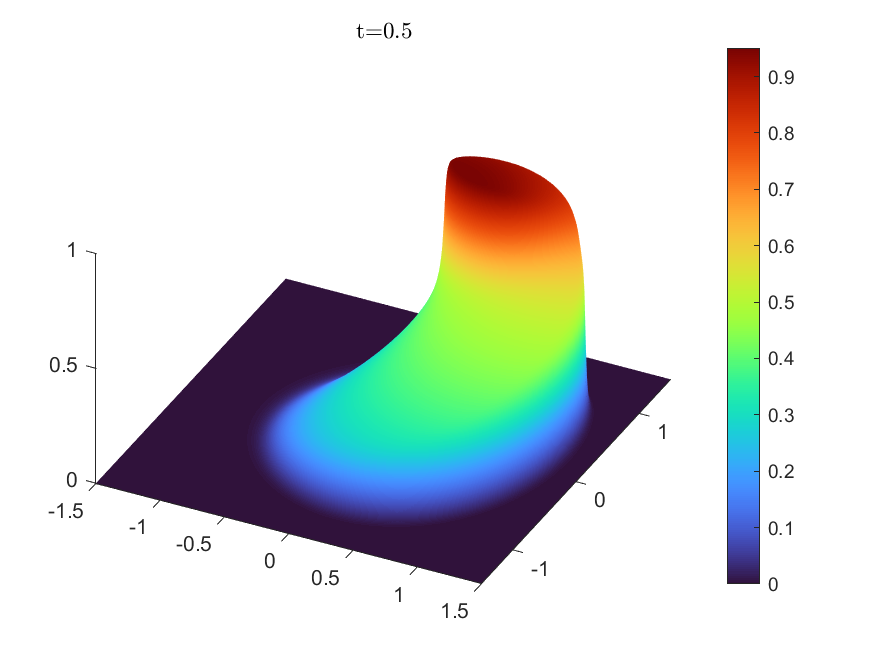} 
    \end{center} 
    \end{minipage}
    \caption{Numerical approximation of the Buckley-Leverett with diffusion \eqref{eq_BL_2D} at final time $t=0.5$.}
    \label{Test4}
\end{figure}

\section{Conclusion}
\label{Section::Conclusion}
In this paper, we have proposed an explicit scheme for the solution of discrete kinetic approximations of advection-diffusion equations. Indeed, the stiffness due to the presence of a relaxation operator could be treated with a projective integration approach, in order to avoid implicit time integration.
The idea of using projective integration methods for discrete kinetic approximation could be easily extended to other types of equations. Indeed, relaxation schemes have a deep connection with lattice Boltzmann Methods (LBM), see for instance \cite{graille2014, wissocqabgrall2024, aregba2024, aregbabellotti2025}. In a very recent paper, the theory of discrete kinetic BGK approximation of hyperbolic problems has been used to prove monotonicity properties for the LBM with overrelaxation \cite{aregba2024}. This approach could be extended to the case of possibly degenerate parabolic equations and the computational performance of LBM could be compared (or coupled) with the projective integration approach. \\

\subsection*{Acknowledgement}
The author would like to thank Thomas Rey and Giuseppe Visconti for comments on this manuscript and for their careful reading. He would also like to thank Roberto Natalini and  Giovanni Samaey for fruitful discussions on this topic.\\ 
The author is funded by the European Union’s Horizon Europe research and innovation programme under the Marie Skłodowska-Curie Doctoral Network DataHyking (Grant No.101072546) and he is member of GNCS-INdAM research group.

\subsection*{Data availability}
Data sharing not applicable to this article as no datasets were generated or analysed during the current study.

\section*{Declarations}
\subsection*{Conflict of interest}
The author declares that he has no financial or non-financial conflicts of interest that are relevant to the content of this article.

\bibliography{References}
\bibliographystyle{acm} 
\appendix
\section{Runge-Kutta Methods}
\label{Appendix_RK}
In Section \ref{Section::Projective} we have presented a higher order extension of projective integration schemes based on Runge-Kutta (RK) methods. Let us recall the main properties of RK methods. We consider an explicit $S$-stage RK method for a general equation of the form
\begin{equation}
    y'=f(t,y), \qquad y(t_0)=y_0.
\end{equation}
Let us consider a time step discretization given by $\Delta t$, then  we have
\begin{equation}
    \begin{cases}
        f^{n+c_s} = f^n + \Delta t \sum_{l=1}^{s-1} a_{s,l} \mathbf{k}_l,\\
        \mathbf{k}_s = f \left(t^n+c_s\,\Delta t, f^{n+c_s} \right)
    \end{cases} \quad 1\leq s \leq S.
\end{equation}
The matrix $\mathbf{A}=(a_{s,l})_{s,l=1}^S$ is the so-called \textit{Runge-Kutta matrix}, the vector $\mathbf{b}=(b_s)_{s=1}^S$ and $\mathbf{c}=(c_s)_{s=1}^S$ are the \textit{Runge-Kutta weights} and the \textit{Runge-Kutta nodes}, respectively. The weights $\mathbf{b}_s$ and $\mathbf{c}_s$ are chosen to ensure consistency of the method, requiring $0 \leq b_s \leq 1$ and $0 \leq c_s \leq 1$ and
\begin{equation}
\label{Assumptions_RK}
    \sum_{s=1}^S b_s = 1, \sum_{s=1}^{S-1} a_{s,l} = c_s, \quad 1 \leq s \leq S.
\end{equation}
Here we represent the Bucther's tableaux used in this paper, where we first recall the general form:
\begin{equation*}
    \begin{tabular}
{c|c}
$c$ &
$A$\\[3pt]
\hline
\rule{0pt}{3ex} 
& $b^T$
\end{tabular}
\end{equation*}

Now, we give detailed expressions of the Runge-Kutta tableaux for the third and the fourth order methods.
\begin{itemize}
        \item Third-order SSP Runge-Kutta
        \begin{equation*}
            \begin{tabular}
                {c|ccc}
                $0$\\
                $1$ & $1$\\
                $1/2$ & $1/4$ & $1/4$\\
                \hline
                & $1/6$ & $1/6$ & $2/3$
            \end{tabular}
        \end{equation*} 
        
    \item Classical fourth-order Runge-Kutta method
        \begin{equation*}
            \begin{tabular}
                {c|cccc}
                $0$\\
                $1/2$ & $1/2$\\
                $1/2$ & $0$ & $1/2$\\
                $1$ & $0$ & $0$ & $1$\\
                \hline
                & $1/6$ & $1/3$ & $1/3$ & $1/6$
            \end{tabular}
        \end{equation*} 
    
\end{itemize}

\section{Stability of the inner integrators}
\label{Appendix_Stability}
In this section we prove the result on the spectrum of the inner integrator, partially retracing the proof in \cite{lafittelejon2016}. The main difference are the non-constant $\eps$-dependent coefficient in front of the entries $D_j$ and the values of the Maxwellians $M_j$.

\begin{oneshot}{\ref{thm:Spectrum_Inner}}
The spectrum of the matrix $\mathcal{A}$ defined in Section \ref{spectrum_inner_integrator} satisfies
\begin{equation}
    \text{Sp}(\mathcal{A}) \subset \mathcal{D} \left(0, C\sqrt{\eps}(|\xi|+|\gamma|) \right) \cup \left\{ \zeta(\eps)\right\},
\end{equation}
where $C$ is a constant depending on the entries of $D$, $\xi$ and $\gamma$ are defined in \eqref{D_matrix_spectrum} and $\zeta(\eps)$ is the dominant eigenvalue.
\end{oneshot}

\begin{proof}
    In the first step we use Rouché's Theorem to localize the eigenvalues, whereas in the second step we give the explicit asymptotic expansion for the dominant eigenvalue. We have obtained the following expression for the characteristic polynomial of $\mathcal{A}$,
    \begin{equation}
    \label{characteristic_polynomial_A}
        \chi_\mathcal{A}(\zeta) = \prod_{i=1}^4 (-\eps D_i - \zeta) \, \left(1-\frac{1}{4}\,\sum_{j=1}^4 \frac{M_j}{\zeta-\eps D_j} \right).
    \end{equation}
    The expression for the characteristic polynomial of matrix $\mathcal{A}_0 := D+P$ is given instead by
    \begin{equation}
        \chi_{\mathcal{A}_0}(\zeta) = \prod_{i=1}^4 (-\eps D_i - \zeta) \, \left(1-\frac{1}{4}\,\sum_{j=1}^4 \frac{1}{\zeta-\eps D_j} \right).
    \end{equation}
    The idea is to study the behavior of the rational function 
    \begin{equation}
        \mathcal{F}: \zeta \to \displaystyle \displaystyle \frac{\chi_{\mathcal{A}(\zeta)}-\chi_{\mathcal{A}_0}(\zeta)}{\chi_{\mathcal{A}_0}(\zeta)} = -\frac{\displaystyle \frac{1}{4} \, \displaystyle \sum_{j=1}^4 (M_j-1) \displaystyle \frac{1}{\zeta-\eps D_j}}{1-\displaystyle \frac{1}{4}\, \displaystyle \sum_{j=1}^4 \displaystyle \frac{1}{\zeta-\eps D_j}}
    \end{equation}
    and to observe the regions containing eigenvalues of $\mathcal{A}_0$ for which $|\mathcal{F}(\zeta)| < 1$.\\
    Here, we just point out some aspects of the computation which differ from the original one.\\
    \begin{itemize}
        \item Let us first consider the dominant eigenvalue of $\mathcal{A}_0$, restricting the analysis to the region parametrized by $\Omega_1 := \left\{ 1 + \eps\,r\, e^{i\theta}, \theta \in [0,2\pi) \right\}$, performing a Taylor expansion of $\displaystyle \frac{1}{1 + \eps\,r\, e^{i\theta} - \eps \, D_j}$ in terms of $\eps$:
        \begin{equation*}
            \frac{1}{1 + \eps\,r\, e^{i\theta} - \eps \, D_j} = 1 + \left(\eps\,D_j - \eps \, r\, e^{i\theta} \right) + \mathcal{O}(\eps^2).
        \end{equation*}
        We observe that 
        \begin{equation}
        \label{Rouche_sum_Maxwellians}
            \sum_{j=1}^4 (M_j-1) = 0.
        \end{equation}
        The rational function $\mathcal{F}$ could be expanded as
        \begin{equation}
            \mathcal{F}(1+\eps\,r\,e^{i\,\theta}) = -\frac{\displaystyle \frac{1}{4} \, \displaystyle \sum_{j=1}^4 (M_j-1)\,D_j + \mathcal{O}(\eps)}{r\,e^{i\,\theta}-\displaystyle \frac{1}{4}\, \displaystyle \sum_{j=1}^4 \displaystyle D_j + \mathcal{O}(\eps)}.
        \end{equation}
        Thus, choosing $r > \displaystyle \frac{1}{4} \left| \sum_{j} D_j \right| + \frac{1}{2} \left| \sum_{j} D_j\,(M_j-1) \right|$, we ensure that $\left|\mathcal{F}(\zeta) \right| < 1/2 + \mathcal{O}(\eps)$.
        From Rouché's Theorem, in $\Omega_1$ (around $\zeta = 1$), $\chi_{\mathcal{A}}$ and $\chi_{\mathcal{A}_0}$ have the same number of zeroes, which means that $\mathcal{A}$ have only one eigenvalue in $\Omega_1$.\\

        \item We focus now on the remaining eigenvalues, by considering the region around $\zeta=0$, which contains the ``fast'' eigenvalues of $\mathcal{A}_0$. Then, we restrict the analysis to the region $\Omega_2 := \left\{ \eps r\,e^{i\,\theta}, \, \theta \in [0,2\pi) \right\}$ and we perform a Taylor expansion of $\displaystyle \frac{1}{r\,e^{i\,\theta}-D_j}$ in terms of $1/r$:
        \begin{equation*}
            \frac{1}{r\,e^{i\,\theta}-D_j} = \frac{e^{-i\theta}}{r} + \frac{D_j \, e^{-i\,2\,\theta}}{r^2} + \mathcal{O} \left(\frac{1}{r^3} \right).
        \end{equation*}
        Using again \eqref{Rouche_sum_Maxwellians}, the rational function $\mathcal{F}$ could be now expanded as
        \begin{multline*}
            \mathcal{F}(\eps\,r\,e^{i\,\theta}) = -\frac{\displaystyle \frac{1}{4} \, \displaystyle \sum_{j=1}^4 (M_j-1)\,D_j \frac{e^{-i\,2\,\theta}}{r^2} + \mathcal{O} \left(\frac{1}{r^3}\right) }{4\eps - \displaystyle \frac{e^{-i\,\theta}}{r} + \displaystyle \mathcal{O} \left(\frac{1}{r^2} \right)} = \\ -\frac{e^{-i\,\theta}}{r} \frac{\displaystyle \frac{1}{4} \, \displaystyle \sum_{j=1}^4 (M_j-1)\,D_j + \mathcal{O} \left(\frac{1}{r}\right) }{4\eps\,r\,e^{i\,\theta} - 1 + \displaystyle \mathcal{O} \left(\frac{1}{r} \right)}.
        \end{multline*}
        We finally perform a Taylor expansion in terms of $\eps$, to obtain
        \begin{equation}
            \mathcal{F}(\eps\,r\,e^{i\,\theta}) = \frac{e^{-i\,\theta}}{r} \displaystyle \frac{1}{4} \, \displaystyle \sum_{j=1}^4 (M_j-1)\,D_j + \mathcal{O}(\eps).
        \end{equation}
        Again, choosing $r > \displaystyle \frac{1}{2} \left| \sum_{j} D_j \right| + 1$, we ensure the inequality $|\mathcal{F}(\zeta) | < 1/2$ and we conclude using Rouché's Theorem that in $\Omega_2$ (around $\zeta=0$), $\chi_{\mathcal{A}}$ and $\chi_{\mathcal{A}_0}$ have the same number of zeroes, which means that $\mathcal{A}$ have $J-1$ eigenvalues in $\Omega_2$.
    \end{itemize}

    In the second step we give an asymptotic expression for the largest eigenvalue in terms of $\eps$. We start from the expression of $\chi_\mathcal{A}$ given in \eqref{characteristic_polynomial_A}. If $\zeta(\eps)$ is a root of $\chi_{\mathcal{A}}$, it holds
    \begin{equation}
    \label{chi_tilde_A}
        \tilde{\chi}_\mathcal{A} (\zeta) := 1-\frac{1}{4} \sum_{j=1}^4 \frac{M_j}{\zeta - \eps \, D_j} = 0.
    \end{equation}
    Let us consider an expression of $\zeta$ in terms of its real coordinates, namely
    \begin{equation}
        \zeta (\eps) = x(\eps) + i\,y(\eps)
    \end{equation}
    and let us consider the asymptotic expansion of $x(\eps)$ an $y(\eps)$ in terms of $\eps$:
    \begin{equation*}
        x(\eps)=x_0+\eps\,x_1+\eps^2\,x_2 + \mathcal{O}(\eps^3), \qquad y(\eps)=y_0+\eps\,y_1+\eps^2\,y_2 + \mathcal{O}(\eps^3).
    \end{equation*}
    We plug these expansions in \eqref{chi_tilde_A} and we perform a Taylor series expansion in terms of $\eps$. Since $\zeta(\eps)$ is a root for $\tilde{\chi}_\mathcal{A}$, all the coefficients of equal power in $\eps$ must vanish. For the zeroth order term we have
    \begin{equation}
        1-\frac{1}{4} \sum_{j=1}^4 \frac{M_j}{x_0+i\,y_0}=0,
    \end{equation}
    which, since $\sum_j (M_j-1)=0$ implies
    \begin{equation*}
        \frac{x_0}{x_0^2+y_0^2}=1 \qquad - \frac{y_0}{x_0^2+y_0^2}=0,
    \end{equation*}
    namely $x_0=1$ and $y_0=0$. Let us now determine the terms of order $\eps$, using
    \begin{equation*}
        \frac{1}{4} \sum_{j=1}^4 M_j \left[(x_1+i\,y_1) + D_j\right]=0,
    \end{equation*}
    which implies   
    \begin{equation*}
        x_1 = \left(1-\frac{1}{\theta^2} \right) \alpha + \frac{1}{\theta^2 \sqrt{\eps}} \xi, \qquad y_1= \frac{1}{\lambda} \beta.
    \end{equation*}
    Finally, let us consider the second order terms in $\eps^2$
    \begin{equation*}
        \frac{1}{4} \sum_{j=1}^4 M_j \left[(x_1+i\,y_1)^2 + 2\,D_j^2 + 4(x_1+i\,y_1)\,D_j-2(x_2+i\,y_2)\right]=0,
    \end{equation*}
    for which we obtain
    \begin{equation*}
        x_2 = 2(x_1^2-y_1^2) + 2\left(1-\frac{1}{\theta^2}\right)(\alpha^2-\beta^2) + \frac{2}{\eps \theta^2} (\xi^2-\gamma^2),
    \end{equation*}
    \begin{equation*}
        y_2 = 4x_1\,y_1 + 2\frac{1}{\lambda}\alpha\,\beta,
    \end{equation*}
\end{proof}

\section{Stability of higher order projective integration}
\label{stability_PRK_appendix}
Let us recall here the stability properties of the projective Runge-Kutta methods. In order to compute the linear stability regions, the Dahlquist equation is introduced
\begin{equation}
\label{Dahlquist_linear}
    \dot{y} = \lambda \, y, \qquad \text{Re}(\lambda) < 0,
\end{equation}
whose general solution does not blow up in time. Introducing a general inner integrator for this equation
\begin{equation}
    y^{k+1} = \tau(\lambda \delta t) y^k,
\end{equation}
we call $\tau(\lambda \delta t)$ the \textit{amplification factor} of the inner integrator. Applying explicit time integration to equation \eqref{Dahlquist_linear}, we observe that 
\begin{equation}
    y^{k+1} = \tau^{k+1}\,y_0,
\end{equation}
which means that the inner integrator is linearly stable if $|\tau| \leq 1$. In the case of forward Euler integrator, we have $\tau(\lambda \delta t) = 1 + \lambda\,\delta t$. The goal of this section is to determine the values of $\tau$ for which the projective integration method is also stable. Let us first focus on projective Forward-Euler method applied to \eqref{Dahlquist_linear}, which can be written as
\begin{equation}
    y^{n+1} = \sigma(\tau, \Delta t, \delta t, K) \, y^n,
\end{equation}
where the amplification factor $\sigma$ is given by
\begin{equation}
\label{sigma_PFE}
    \sigma^{\text{PFE}}(\tau, \Delta t, \delta t, K) = \left[ \left( \frac{\Delta t - (K+1)\,\delta t}{\delta t} + 1 \right) \tau - \frac{\Delta t - (K+1) \delta t}{ \delta t} \right] \tau^k.
\end{equation}
The PFE method is stable if $|\sigma^{\text{PFE}}(\tau, \Delta t, \delta t, K)| \leq 1$. The main goal of introducing a projective integration approach is to use a time step $\Delta t = \mathcal{O}(\Delta x^\gamma)$. In the limit $\eps \to 0$ for fixed $\Delta x$, we look at the limiting stability regions as $\Delta t/ \delta t \to \infty$. In \cite{gear2003} it has been shown that the values $\tau$ satisfying \eqref{sigma_PFE} lie in the union of two separated disks $\mathcal{D}_1^{\text{PFE}} \cup \mathcal{D}_2^{\text{PFE}}$, where
\begin{equation}
    \mathcal{D}_1^{\text{PFE}}= \mathcal{D} \left(1-\frac{\delta t}{\Delta t}, \frac{\delta t}{\Delta t} \right), \qquad \mathcal{D}_2^{\text{PFE}} \left(0, \left(\frac{\delta t}{\Delta t} \right)^{\frac{1}{K}} \right),
\end{equation}
and $\mathcal{D}(c, r)$ denotes the ball with center $(c,0)$ and radius $r$ in the complex $\lambda$-plane.\\
For higher order projective Runge-Kutta (PRK) schemes, it is possible to give an explicit expression of the \textit{amplification factor}, applying the integrator to Equation \eqref{Dahlquist_linear}. For the sake of simplicity, we introduce 
\begin{equation}
    M=\frac{\Delta t - (K+1) \delta t}{\delta t}, \quad \text{ and } \quad M_s=\frac{c_s\,\Delta t - (K+1) \delta t}{\delta t},
\end{equation}
where $c_s$ are the RK nodes. We also denote by $\mathbf{A}=(a_{s,l})_{s=1}^S$ and $\mathbf{b}=(b_s)_{s=1}^S$ the RK matrix and the RK weights for a general $S$-stage Runge-Kutta method. Then, for an explicit PRK method, we have
\begin{equation}
    \begin{cases}
        \mathbf{k}_1 = \kappa_1(\tau) y^n = \displaystyle \frac{\tau^{K+1}-\tau^K}{\delta t} \, y^n, \\
        \mathbf{k}_s = \kappa_s(\tau) y^n = \displaystyle \frac{\tau^{K+1}-\tau^K}{\delta t} \left(\tau^{K+1} + (M_s\,\delta t) \sum_{l=1}^{s-1} \frac{a_{s,l}}{c_s} \kappa_l \right) \, y^n, \qquad 2 \leq s \leq S,\\
        y^{n+1} = \displaystyle \sigma(\tau)\,y^n = \left( \tau^{K+1} + (M\,\delta t) \sum_{s=1}^S b_s\,\kappa_s \right) \, y^n.
    \end{cases}
\end{equation}
The following result on the stability of higher order projective Runge-Kutta methods holds.
\begin{thm}[\cite{lafittelejon2016}, Theorem 3.2]
\label{Thm:PRK_stability} 
Assume the inner integrator is stable, i.e. $|\tau| \leq 1$ and $\delta t$, $\Delta t$ and $K$ are chosen such that the projective forward Euler method is stable. Then, a projective Runge-Kutta method is also stable if it satisfies the following assumption
\begin{equation}
    \sum_{s=1}^S b_s = 1, \sum_{s=1}^{S-1} a_{s,l} = c_s, \quad 1 \leq s \leq S,
\end{equation}
with $0 \leq b_s, \, c_s \leq 1$ and the convexity condition 
\begin{equation}
    0 \leq a_{s,l} \leq c_s, \qquad \forall\, 1 \leq l \leq s, \quad \forall\, 1 \leq s \leq S.
\end{equation}
\end{thm}

\begin{figure}[h!]
\centering
\subfigure[]{\includegraphics[scale=0.32]{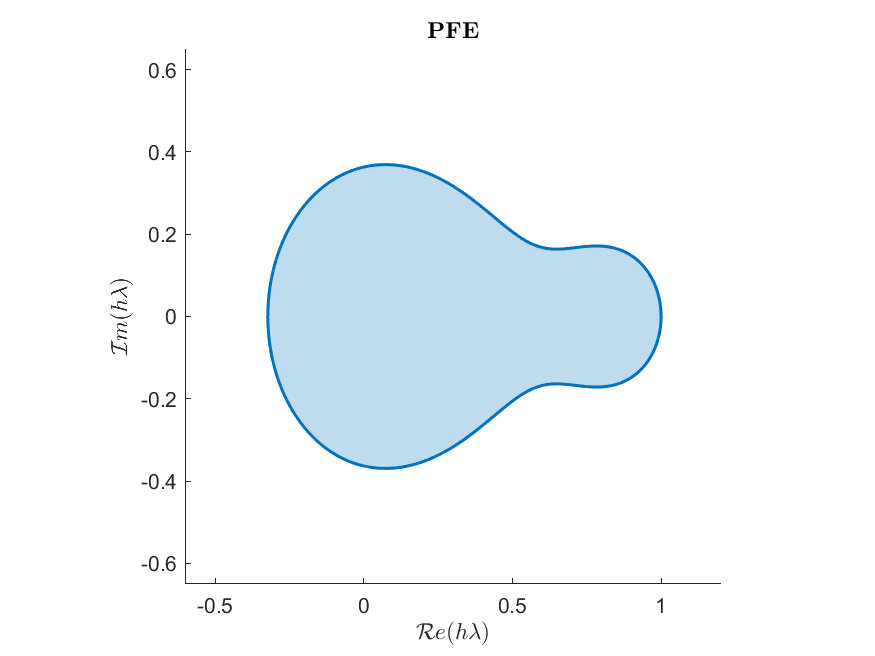}}
\subfigure[]{\includegraphics[scale=0.32]{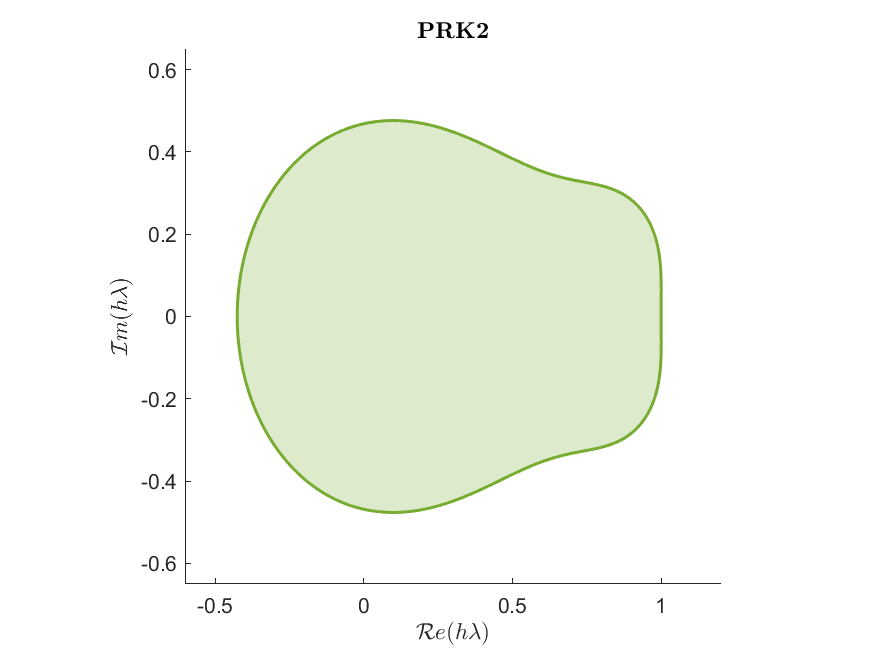}}
\subfigure[]{\includegraphics[scale=0.32]{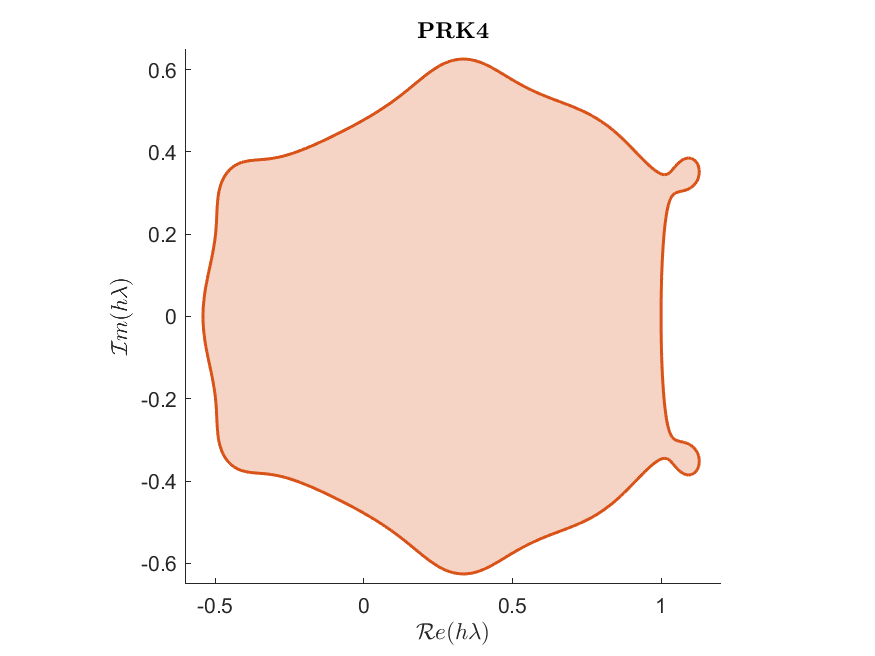}}
\caption{Stability region for higher order projective integrators in the complex $h\lambda$-plane for $\Delta t/\delta t = 10$.}
\label{Fig:PRK_Stability_10}
\end{figure}

\begin{figure}[h!]
\centering
\subfigure[]{\includegraphics[scale=0.31]{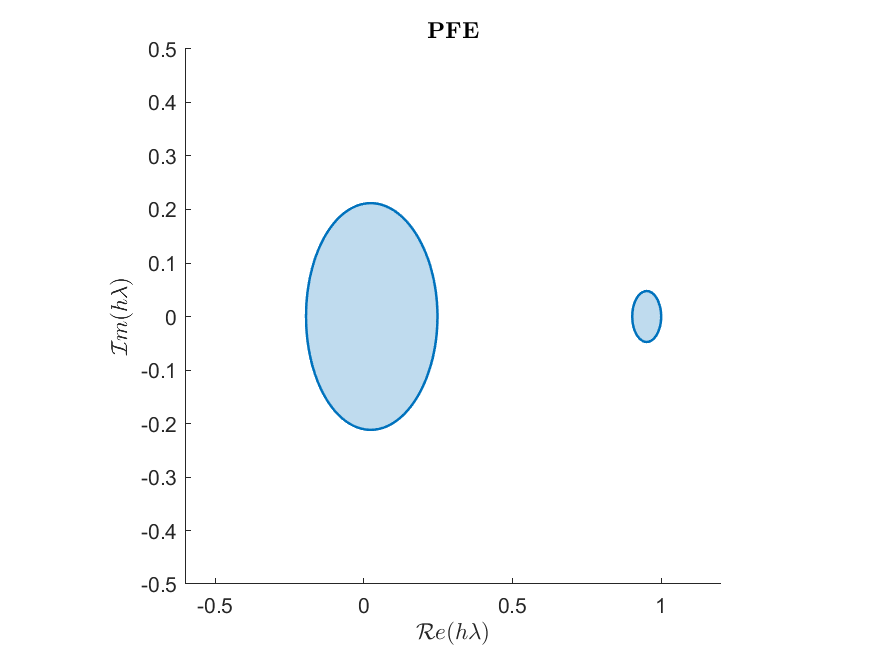}}
\subfigure[]{\includegraphics[scale=0.31]{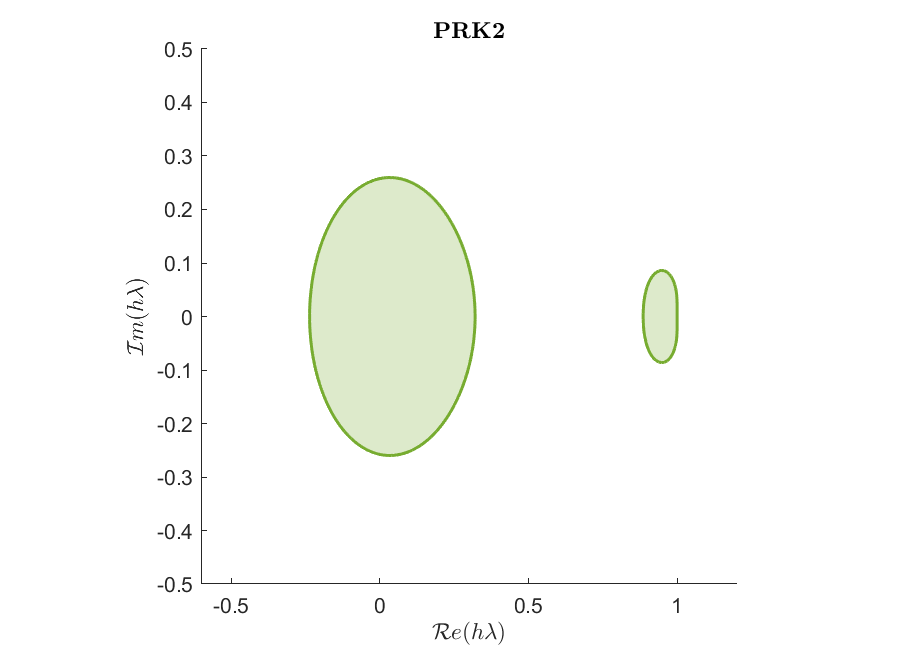}}
\subfigure[]{\includegraphics[scale=0.31]{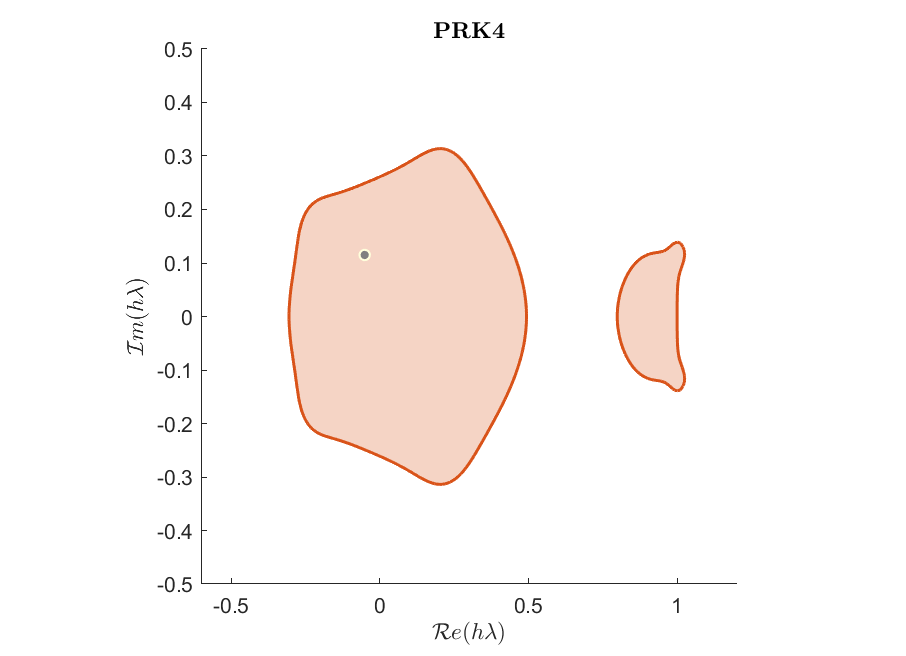}}
\caption{Stability region for higher order projective integrators in the complex $h\lambda$-plane for $\Delta t/\delta t = 25$.}
\label{Fig:PRK_Stability_25}
\end{figure}

In Fig. \ref{Fig:PRK_Stability_10}-\ref{Fig:PRK_Stability_25} we show the stability regions of projective Runge-Kutta schemes for different values of $\delta t/\Delta t$ and $K$. The main feature, as proven in Theorem \ref{Thm:PRK_stability}, is that the qualitative dependence of higher order PRK methods on the parameters is equivalent to those of PFE.\\

\end{document}